\documentclass{article}
\usepackage{authblk}
\usepackage{float}
\usepackage{amsmath,amssymb,amsthm}
\usepackage{graphicx}
\usepackage{array}
\usepackage{booktabs}
\usepackage{longtable}
\usepackage{color}
\usepackage{caption}
\usepackage[super]{natbib}
\newcommand{\fig}[1]{Figure~\ref{#1}}

\newcommand{\tab}[1]{Table~\ref{#1}}

\title{\bf Monitoring the Ratio of two Normal Variables using EWMA Type Control Charts in Short Production Runs}
\author[2]{Thi Hien Nguyen}
\author[3]{Jean-Michel Masereel}
\author[4]{Guillaume Tartare}
\author[1,4]{Kim Duc Tran\thanks{Corresponding author: Kim Duc Tran. Email: \texttt{ductk@donga.edu.vn}}}

\affil[1]{International Chair in DS \& XAI, International Research Institute for Artificial Intelligence and Data Science, Dong A University, Danang, Vietnam}

\affil[2]{Laboratoire AGM, UMR CNRS 8088, CY Cergy Paris Université, 95000 Cergy, France}

\affil[3]{ESIEE-IT, École d'Ingénieurs et d'Experts IT Informatique, 95300 Pontoise, France}

\affil[4]{Univ. Lille, ENSAIT, ULR 2461 -- GEMTEX -- Génie et Matériaux Textiles, F-59000 Lille, France}

\begin{document}

%%%%%%%%%%%%%%%%%%%%%%%%%%%%%%%%%%%%%%%%%%%%%%%%%%%%%%%%%%%%%%%%%%%%%%%%%%%%%%%%
\maketitle

\begin{abstract}
In many industrial and engineering applications, process performance is characterized by the ratio of two normally distributed quality characteristics. Monitoring such ratios is particularly challenging in short production cycles, where conventional control charts often suffer from limited sensitivity due to the small number of available inspections. This paper proposes an Exponentially Weighted Moving Average (EWMA) type control chart for monitoring the ratio of two normally distributed random variables under short production run (SPR) conditions. The statistical distribution of the ratio is first reviewed, adopting the corrected closed-form density of \citet{Nadarajah2020} rather than the approximation used in earlier studies. The control limit of the proposed chart is calibrated to a prescribed in-control truncated average run length (TARL$_0$) over a finite horizon $I$ of inspections, using a Markov-chain representation of the EWMA recursion. The detection performance of the chart is then assessed through a large factorial study covering the smoothing constant $\lambda$, the in-control correlation $\rho_0$, the coefficients of variation $(\gamma_X,\gamma_Y)$, the sample size $n$, and the magnitude of the shift $\tau$. Numerical results show that the proposed EWMA-RZ chart provides substantially better detection of small and moderate shifts than the recently developed Shewhart-type short-run chart for the ratio (ShRZ) of \citet{Tran2021}, especially for $|\tau-1|\le 0.05$. An illustrative example based on a beverage filling process is included to demonstrate the practical implementation of the method.
\end{abstract}


%%%%%%%%%%%%%%%%%%%%%%%%%%%%%%%%%%%%%%%%%%%%%%%%%%%%%%%%%%%%%%%%%%%%%%%%%%%%%%%%
\section{Introduction}
\label{sec:introduction}

Following the major events humanity has recently experienced, from the global COVID-19 pandemic to severe natural disasters, the urgency of adopting ``smart manufacturing'' and clear quality management processes has become even more apparent. For example, the 7.7-magnitude earthquake that struck Sagaing, Mandalay, Myanmar, in March 2025 caused extensive damage, resulting in a large number of deaths and injuries. Simultaneously, the final months of 2025 saw a series of severe floods and landslides in several South Asian and Southeast Asian countries (including Indonesia, Thailand, Sri Lanka and Vietnam), killing more than 1{,}100 people and causing significant damage to infrastructure and transportation. The combined impact of pandemics and natural disasters, not only on human lives but also on the economy and society, creates an urgent need for robust production, healthcare and relief systems. Effective and stable quality control processes are required to minimize damage and enhance response capabilities. In this context, Statistical Process Control (SPC), through traditional statistical tools or in combination with modern methods such as Artificial Intelligence, Machine Learning, and Deep Learning, enables the monitoring, analysis, and control of abnormal fluctuations during processes. One of the most practical and commonly used applications is the use of control charts, which help detect deviations early and ensure process stability and reliability.

Sustainable products have always held a significant place in the modern market. However, in the context of digitalization and increasingly diverse and flexible customer demands, traditional production models struggle to keep pace with the rate of change. This has led to the emergence and development of the \emph{short production run} (SPR) model. Economically, a short run is understood as a phase in which at least one factor of production, usually capital, is kept constant, while other factors, such as labor, can change. In the manufacturing sector, a short production run refers to product batches with small quantities, short production times, and highly flexible adjustments to product design or configuration. Unlike mass production, this model focuses on responding quickly to fluctuations in market demand, testing new products, or fulfilling customized orders \citep{Tang1994}. By reducing machine setup time, optimizing processes, and increasing flexibility, short production runs help businesses shorten product development cycles, reduce inventory, and mitigate risks arising from constantly changing demand \citep{Celano2011,Li2012}. Several studies have shown that SPR reduces costs associated with inventory, obsolescence, and forecasting errors, particularly in industries with short product life cycles \citep{Celano2011,Khoo2022,Tran2021}. In this context, contributions include EWMA charting \citep{Ong2023}, adaptive charts for varying sample sizes \citep{Castagliola2013,Amdouni2015}, and economic designs \citep{Celano2012,Zhang2014}. Monitoring specific parameters such as coefficients of variation \citep{Amdouni2016,Amdouni2017,Khatun2019,Chew2025} or ratios of normal variables \citep{Tran2021} has also been the subject of recent work. These studies highlight the importance of developing robust statistical tools capable of rapidly detecting deviations under data-scarce conditions \citep{Khoo2022,Sfiris2014}. \tab{tab:recent_work} summarizes recent control-chart methods designed for short production runs.

\begin{table}[!htbp]
	\centering
	\caption{Summary of control-chart methods for short production runs}\label{tab:recent_work}
	\begin{tabular}{p{0.25\linewidth}p{0.27\linewidth}p{0.40\linewidth}}
		\toprule
		\textbf{Reference} & \textbf{Chart type} & \textbf{Parameter / Method} \\
		\midrule
		\citet{Celano2011}      & Shewhart, EWMA         & $t$-statistic (mean) \\
		\citet{Li2012}          & Various                & Two-sided control charts \\
		\citet{Castagliola2013} & VSS $t$-chart          & Adaptive sampling for mean \\
		\citet{Celano2012}      & CUSUM $t$-chart        & Economic performance analysis \\
		\citet{Zhang2014}       & SPRT                   & Economically designed chart \\
		\citet{Amdouni2015}     & VSS                    & Coefficient of variation \\
		\citet{Amdouni2016}     & Run rules              & One-sided CV monitoring \\
		\citet{Sfiris2014}      & Adaptive fuzzy         & Fuzzy estimators \\
		\citet{Amdouni2017}     & VSI Shewhart           & CV with variable intervals \\
		\citet{Chong2019}       & Hotelling's $T^2$      & Multivariate with VSS \\
		\citet{Khatun2019}      & One-sided              & Multivariate CV \\
		\citet{Khoo2022}        & Median chart           & Robust location estimator \\
		\citet{Ong2023}         & EWMA for $S^2$         & Variance monitoring \\
		\citet{Tran2021}        & One-sided Shewhart     & Ratio of two normals \\
		\citet{Chew2025}        & VSS                    & Multivariate CV optimization \\
		\bottomrule
	\end{tabular}
\end{table}

It is important to distinguish short production runs from seasonal products. While short production runs may have some negative environmental impacts, such as products being easily discarded and the high frequency of machine setups leading to increased energy consumption and waste, this model also offers significant benefits. Short production runs ensure that production meets demand, reduce inventory levels, and limit product obsolescence. If businesses adopt lean manufacturing practices, integrate sustainable production planning with SPC, and implement sustainable design based on a circular economy, short production runs can even become a more environmentally friendly option than traditional mass production.

The statistical monitoring of the ratio of two normally distributed variables has gained significant attention in recent quality-control research due to its practical relevance in various applications. Initial developments by \citet{Tran2016} established EWMA-based control charts specifically designed for this parameter, addressing the non-normality challenges inherent in ratio distributions. Subsequent research introduced adaptive sampling strategies, with \citet{Nguyen2019} proposing variable sampling interval EWMA charts and \citet{Nguyen2020} developing CUSUM charts with similar adaptive features, both demonstrating improved performance over fixed sampling schemes. The latest advancements by \citet{Haq2024} focus on enhanced memory-type control charts that offer superior detection capabilities for shifts in the ratio parameter. These methodologies provide quality practitioners with sophisticated tools for monitoring critical ratios in processes where the relationship between two variables is more informative than their individual behaviors, such as in yield monitoring, cost-efficiency ratios, or technical performance indicators.

To the best of our knowledge, however, no contribution has yet combined the memory-type sensitivity of the EWMA chart with the finite-horizon nature of short production runs for the ratio of two normal variables. The closest references, \citet{Tran2021} and \citet{Khoo2022}, develop Shewhart-type short-run charts; while their constructions are tractable, their sensitivity to small and moderate shifts is limited, which is precisely the regime in which the EWMA chart is known to outperform Shewhart-type schemes \citep{Tran2016,Nguyen2019,Haq2024}. The present paper fills this gap by proposing an EWMA-type chart for the ratio $Z=X/Y$ designed for SPR conditions, using the truncated average run length (TARL) as the design and assessment criterion.

The remainder of this paper is organized as follows. Section~\ref{sec:distribution} presents the distribution of the ratio $Z$ between two normally distributed random variables. We adopt the corrected closed-form expression of \citet{Nadarajah2020} for the density of $Z$, rather than the approximation used in \citet{Castagliola2013} and \citet{Nguyen2019}. Section~\ref{sec:implementation} formulates the proposed EWMA-RZ chart for short production runs. Section~\ref{sec:optimization} describes the Markov-chain-based computation of the TARL and the design procedure for calibrating the upper control limit. Section~\ref{sec:numerical} reports the numerical results of an extensive factorial study and compares the EWMA-RZ chart with the ShRZ chart of \citet{Tran2021}. An illustrative example is given in Section~\ref{sec:illustrative}, and concluding remarks are presented in Section~\ref{sec:conclusions}.

%%%%%%%%%%%%%%%%%%%%%%%%%%%%%%%%%%%%%%%%%%%%%%%%%%%%%%%%%%%%%%%%%%%%%%%%%%%%%%%%
\section{The distribution of the ratio $Z$}
\label{sec:distribution}

Suppose that $X$ and $Y$ are two normal random variables such that $\mathbf{W}=(X,Y)^T \sim N(\pmb{\mu}_\mathbf{W},\pmb{\Sigma}_\mathbf{W})$, i.e.\ $\mathbf{W}$ is a bivariate normal random vector with mean vector and variance-covariance matrix
\begin{align}\label{eq:bvn}
	{\pmb{\mu}}_{\mathbf{W}}=\begin{pmatrix}\mu_{X}\\ \mu_{Y}\end{pmatrix},
	\qquad
	{\pmb{\Sigma}}_{\mathbf{W}}=\begin{pmatrix}
		\sigma^2_{X}            & \rho\,\sigma_{X}\sigma_{Y}\\
		\rho\,\sigma_{X}\sigma_{Y} & \sigma^2_{Y}
	\end{pmatrix},
\end{align}
where $\mu_X$ and $\mu_Y$ are the means of the two variables and $\rho$ is their correlation coefficient. The coefficients of variation $(\gamma_X,\gamma_Y)$ and the standard-deviation ratio $\omega$ are defined as
\begin{align*}
	\gamma_X = \frac{\sigma_X}{\mu_X},
	\qquad
	\gamma_Y = \frac{\sigma_Y}{\mu_Y},
	\qquad
	\omega   = \frac{\sigma_X}{\sigma_Y}.
\end{align*}

Let $Z=X/Y$ be the ratio of $X$ to $Y$. \citet{Nadarajah2020} derived a closed-form expression for the probability density function of $Z$ which corrects earlier formulations and avoids the use of approximations to the cumulative distribution function (c.d.f.) of $Z$ as a function of $(\gamma_X,\gamma_Y,\omega,\rho)$ such as the one used in \citet{Celano2016_Synthentic_RZ}. The exact density reads
\begin{equation}
	\label{eq:fz}
	\begin{aligned}
		f_Z(z)
		&=
		\frac{1}{2\pi\,\sigma_X\,\sigma_Y\,\sqrt{1-\rho^2}}\;
		\exp\!\left(
		-\frac{B^2}{2\sigma_X^2(1-\rho^2)}
		-\frac{\mu_Y^2}{2\sigma_Y^2}
		\right)\\[0.4em]
		&\quad\times
		\left[
		\frac{1}{\alpha}
		+\frac{\sqrt{\pi}\,\beta}{2\alpha^{3/2}}\,
		\exp\!\left(\frac{\beta^2}{4\alpha}\right)
		\operatorname{erf}\!\left(\frac{\beta}{2\sqrt{\alpha}}\right)
		\right],
	\end{aligned}
\end{equation}
where $\operatorname{erf}(\cdot)$ denotes the error function
\begin{equation}
	\operatorname{erf}(x) = \frac{2}{\sqrt{\pi}}\int_{0}^{x} e^{-t^2}\,dt,
\end{equation}
and the auxiliary quantities are defined by
\begin{align*}
	&A = z - \frac{\rho\,\sigma_X}{\sigma_Y},
	&&B = \mu_X - \frac{\rho\,\sigma_X\,\mu_Y}{\sigma_Y},\\
	&\alpha = \frac{A^2}{2\sigma_X^2(1-\rho^2)} + \frac{1}{2\sigma_Y^2},
	&&\beta = -\frac{A B}{\sigma_X^2(1-\rho^2)} - \frac{\mu_Y}{\sigma_Y^2}.
\end{align*}

For computational convenience in the Markov-chain construction of Section~\ref{sec:optimization}, the c.d.f.\ of $Z$ can either be obtained by numerical integration of \eqref{eq:fz}, or evaluated through the closed-form approximation
\begin{align}\label{eq:Zcdf}
	F_Z(z\mid\gamma_X,\gamma_Y,\omega,\rho)
	\simeq \Phi\!\left(\frac{A^\star}{B^\star}\right),
\end{align}
with
\begin{align*}
	A^\star = \frac{z}{\gamma_Y} - \frac{\omega}{\gamma_X},
	\qquad
	B^\star = \sqrt{\omega^2 - 2\rho\omega z + z^2},
\end{align*}
where $\Phi(\cdot)$ is the c.d.f.\ of the standard normal distribution. The approximation in \eqref{eq:Zcdf} is the one used in \citet{Celano2016_Synthentic_RZ} and \citet{Nguyen2019}; we use it here as a fast surrogate for the Markov-chain build, while \eqref{eq:fz} is retained as the reference expression for high-accuracy validation. The corresponding approximate inverse distribution function (i.d.f.) is
\begin{align}\label{eq:Zidf}
	F^{-1}_{Z}(p\mid\gamma_X,\gamma_Y,\omega,\rho)
	\simeq
	\begin{cases}
		\dfrac{-C_2-\sqrt{C_2^2-4 C_1 C_3}}{2 C_1}, & p\in(0,0.5],\\[0.5em]
		\dfrac{-C_2+\sqrt{C_2^2-4 C_1 C_3}}{2 C_1}, & p\in[0.5,1),
	\end{cases}
\end{align}
with
\begin{align*}
	C_1 &= \frac{1}{\gamma_Y^2}-\bigl(\Phi^{-1}(p)\bigr)^2,\\
	C_2 &= 2\omega\!\left(\rho\bigl(\Phi^{-1}(p)\bigr)^2 - \frac{1}{\gamma_X\gamma_Y}\right),\\
	C_3 &= \omega^2\!\left(\frac{1}{\gamma_X^2}-\bigl(\Phi^{-1}(p)\bigr)^2\right),
\end{align*}
and $\Phi^{-1}(\cdot)$ is the standard normal quantile function.

%%%%%%%%%%%%%%%%%%%%%%%%%%%%%%%%%%%%%%%%%%%%%%%%%%%%%%%%%%%%%%%%%%%%%%%%%%%%%%%%
%%  >>>>> INSERTION 1: paragraph at end of Section 2 <<<<<
%%%%%%%%%%%%%%%%%%%%%%%%%%%%%%%%%%%%%%%%%%%%%%%%%%%%%%%%%%%%%%%%%%%%%%%%%%%%%%%%
\paragraph{On the bivariate normal assumption.}
The bivariate normal model for $(X,Y)$ underlying the closed-form ratio density~\eqref{eq:fz} is the standard modelling convention in the literature on ratio-of-two-normals control charts \citep{Celano2016_Synthentic_RZ,Tran2016,Nguyen2019,Nguyen2020,Tran2021,Haq2024}, and it is what makes the exact density of \citet{Nadarajah2020} and the closed-form approximation~\eqref{eq:Zcdf} available. In practice this assumption is only ever approximately valid: real quality characteristics may exhibit mild skewness or heavier tails than the normal model permits, particularly when the coefficients of variation are not negligible. The robustness of the proposed EWMA-RZ$^+$ chart to such departures from bivariate normality is examined empirically in Section~\ref{ssec:robustness}.

%%%%%%%%%%%%%%%%%%%%%%%%%%%%%%%%%%%%%%%%%%%%%%%%%%%%%%%%%%%%%%%%%%%%%%%%%%%%%%%%
\section{The EWMA-RZ chart for short production runs}
\label{sec:implementation}

The production run is planned to produce a small-size lot consisting of $N$ parts after a fixed rolling length $H$. Let $I$ be the number of planned inspections over the rolling horizon $H$, and assume that no inspection takes place at the end of the run. Under these settings, the sampling frequency (the time interval between two consecutive inspections) is $\mathcal{S}_h = H/(I+1)$ hours. At each inspection $i=1,2,\dots,I$, a sample of $n$ items is collected and the quality characteristic $\mathbf{W}$ is measured for each item. Let $[\mathbf{W}_{i,1},\dots,\mathbf{W}_{i,n}]$ denote the collected sample, where the couples $\mathbf{W}_{i,j}=(X_{i,j},Y_{i,j})^T$, $j=1,\dots,n$, follow the bivariate normal model $N(\pmb{\mu}_{\mathbf{W},i},\pmb{\Sigma}_{\mathbf{W},i})$ with
\begin{align}\label{eq:bvn_i}
	{\pmb{\mu}}_{\mathbf{W},i}=\begin{pmatrix}\mu_{X,i}\\ \mu_{Y,i}\end{pmatrix},
	\qquad
	{\pmb{\Sigma}}_{\mathbf{W},i}=\begin{pmatrix}
		\sigma^2_{X,i}              & \rho\,\sigma_{X,i}\sigma_{Y,i}\\
		\rho\,\sigma_{X,i}\sigma_{Y,i} & \sigma^2_{Y,i}
	\end{pmatrix},
	\quad i=1,2,3,\dots
\end{align}

Let $\gamma_X$ and $\gamma_Y$ denote the known and constant coefficients of variation, and let $z_0=\mu_{X,i}/\mu_{Y,i}$ and $\rho_0$ be the known in-control values of the ratio and of the correlation coefficient, respectively, that ensure the stability of the process. We further assume the linear relations $\sigma_{X,i}=\gamma_X\,\mu_{X,i}$ and $\sigma_{Y,i}=\gamma_Y\,\mu_{Y,i}$ for every $i\ge 1$, which imply that the standard deviation of each sample changes proportionally to its mean so that their ratio remains constant. This assumption is appropriate for several quality characteristics encountered in practice, such as weights, tensile strengths, and linear dimensions, whose dispersion is proportional to the population mean. The sample units are otherwise free to vary across samples; in particular, $\pmb{\mu}_{\mathbf{W},i}\neq\pmb{\mu}_{\mathbf{W},k}$ and $\pmb{\Sigma}_{\mathbf{W},i}\neq\pmb{\Sigma}_{\mathbf{W},k}$ for $i\neq k$ are allowed. The monitoring statistic of the proposed chart is the ratio of sample means
\begin{align}\label{eq:Zhat}
	\hat Z_i = \frac{\bar X_i}{\bar Y_i}
	= \frac{\sum_{j=1}^{n} X_{i,j}}{\sum_{j=1}^{n} Y_{i,j}},
	\qquad i=1,2,\dots,I.
\end{align}

To obtain the c.d.f.\ and i.d.f.\ of $\hat Z_i$, recall that $\bar X_i \sim N(\mu_{X,i},\sigma_{X,i}/\sqrt{n})$ and $\bar Y_i \sim N(\mu_{Y,i},\sigma_{Y,i}/\sqrt{n})$, with constant coefficients of variation $\gamma_{\bar X}=\gamma_X/\sqrt{n}$ and $\gamma_{\bar Y}=\gamma_Y/\sqrt{n}$. By construction, the standard-deviation ratio at inspection $i$ satisfies
\begin{align}\label{eq:omega0}
	\omega_i = \frac{\sigma_{X,i}}{\sigma_{Y,i}}
	= \frac{\mu_{X,i}}{\mu_{Y,i}}\cdot\frac{\gamma_X}{\gamma_Y}
	= z_0\,\frac{\gamma_X}{\gamma_Y}
	= \omega_0,
\end{align}
where $\omega_0$ is the in-control standard-deviation ratio. Consequently, the c.d.f.\ and i.d.f.\ of $\hat Z_i$ follow directly from those of $Z$ in \eqref{eq:Zcdf}--\eqref{eq:Zidf} as \citep{Nguyen2019}
\begin{align}
	F_{\hat Z_i}(z\mid n,\gamma_X,\gamma_Y,z_0,\rho_0)
	&= F_Z\!\left(z\,\bigg|\,\tfrac{\gamma_X}{\sqrt{n}},\tfrac{\gamma_Y}{\sqrt{n}},\tfrac{z_0\gamma_X}{\gamma_Y},\rho_0\right),\label{eq:Fzhat}\\
	F^{-1}_{\hat Z_i}(p\mid n,\gamma_X,\gamma_Y,z_0,\rho_0)
	&= F^{-1}_Z\!\left(p\,\bigg|\,\tfrac{\gamma_X}{\sqrt{n}},\tfrac{\gamma_Y}{\sqrt{n}},\tfrac{z_0\gamma_X}{\gamma_Y},\rho_0\right).\label{eq:Fzhatinv}
\end{align}

The proposed \emph{EWMA-RZ$^+$} chart for short production runs plots the upper-sided EWMA statistic
\begin{align}\label{eq:EWMA}
	W_i = \lambda\,\hat Z_i + (1-\lambda)\,W_{i-1},
	\qquad i = 1,2,\dots,I,
\end{align}
initialized at $W_0 = z_0$, against a single upper control limit $\mathit{UCL}$. Here $\lambda\in(0,1]$ is the smoothing constant. The chart signals an out-of-control situation as soon as $W_i\ge \mathit{UCL}$, in which case the production phase is interrupted, the root cause is identified, and corrective actions are taken. The detection of upward shifts ($\tau>1$) is the focus of this paper; the lower-sided case is symmetric and follows the same construction with an appropriate lower limit $\mathit{LCL}$. A two-sided chart can be obtained as the conjunction of the upper- and lower-sided charts. We focus on the upper-sided design because in many practical situations (for instance, when $X$ measures a defect-related quantity and $Y$ a normalizing reference), an increase in $Z$ is the alarm of interest.

%%%%%%%%%%%%%%%%%%%%%%%%%%%%%%%%%%%%%%%%%%%%%%%%%%%%%%%%%%%%%%%%%%%%%%%%%%%%%%%%
\section{TARL optimization for the EWMA-RZ control chart}
\label{sec:optimization}

\subsection{The TARL criterion}

In long production runs, the standard performance criterion for a control chart is the average run length (ARL), defined as the expected number of inspections until the chart issues a signal. Under short-run conditions, however, the production is terminated after $I$ inspections, and the ARL is no longer a well-defined quantity since the run length can take only values in $\{1,2,\dots,I,I+1\}$, where $T=I+1$ encodes the event ``no signal during the run''. Following \citet{Castagliola2013} and \citet{Tran2021}, we adopt the \emph{truncated average run length}
\begin{align}\label{eq:TARL}
	\mathrm{TARL}(I) \;=\; \sum_{k=1}^{I} k\,\Pr(T=k)\;+\;(I+1)\,\Pr(T>I),
\end{align}
where the random variable $T$ denotes the number of inspections until the chart signals (or $I+1$ if no signal occurs). The in-control TARL, denoted $\mathrm{TARL}_0$, is computed under the assumption that the process operates with $E[\hat Z_i] = z_0$ for all $i$; the out-of-control TARL, denoted $\mathrm{TARL}_1$, is computed under a shifted mean $E[\hat Z_i] = \tau\,z_0$, where $\tau\neq 1$ is the shift factor.

The design problem for the EWMA-RZ$^+$ chart is to choose the upper control limit $\mathit{UCL}$ so that, for given $\lambda$, $n$, $\gamma_X$, $\gamma_Y$, $\rho_0$, $z_0$, and $I$, the in-control TARL matches a prescribed target value $\mathit{TARL}_0^*$:
\begin{align}\label{eq:UCL_problem}
	\text{find } \mathit{UCL}\text{ such that } \mathrm{TARL}_0(\mathit{UCL}) \;=\; \mathit{TARL}_0^*.
\end{align}
Typical choices are $\mathit{TARL}_0^* = I$, so that the chart signals roughly once over the run when the process is in control, which is consistent with the choice in \citet{Tran2021} for the ShRZ chart.

\subsection{Markov-chain representation of the EWMA recursion}

To evaluate $\mathrm{TARL}_0$ and $\mathrm{TARL}_1$, we use the standard Markov-chain approximation of the continuous-state EWMA process \citep{Brook1972,Lucas1990}. The safe region for the EWMA statistic is $[w_{\min},\,\mathit{UCL})$, where $w_{\min}$ is a lower bound chosen far enough below $z_0$ to capture the bulk of the EWMA distribution. The interval $[w_{\min},\,\mathit{UCL})$ is partitioned into $m$ subintervals of equal width $h = (\mathit{UCL}-w_{\min})/m$, with midpoints $c_1<c_2<\dots<c_m$. An additional absorbing state, denoted $c_{m+1}$, represents the event ``signal''.

Given the current EWMA value $W_{k}=c_i$, the next value is $W_{k+1} = \lambda\,\hat Z_{k+1} + (1-\lambda)\,c_i$. The transition probability from state $i$ to state $j\in\{1,\dots,m\}$ is therefore
\begin{align}\label{eq:Qij}
	Q_{ij}
	&= \Pr\!\left( c_j - \tfrac{h}{2} \le W_{k+1} < c_j + \tfrac{h}{2}\,\Big|\,W_k=c_i\right)\notag\\
	&= F_{\hat Z}\!\left(\tfrac{c_j+\tfrac{h}{2}-(1-\lambda)c_i}{\lambda}\right)
	  - F_{\hat Z}\!\left(\tfrac{c_j-\tfrac{h}{2}-(1-\lambda)c_i}{\lambda}\right),
\end{align}
where $F_{\hat Z}$ denotes the c.d.f.\ of $\hat Z$ given in \eqref{eq:Fzhat}, evaluated under the in-control parameters for $\mathrm{TARL}_0$ and under the shifted parameters $(\tau z_0,\rho_0)$ for $\mathrm{TARL}_1$. The probability of being absorbed (i.e.\ of issuing a signal) from state $i$ is
\begin{align}\label{eq:palarm}
	p_i^{\mathrm{abs}}
	= 1 - F_{\hat Z}\!\left(\tfrac{\mathit{UCL}-(1-\lambda)c_i}{\lambda}\right).
\end{align}

The c.d.f.\ $F_{\hat{Z}_i}$ appearing in the transition probabilities~(14) 
is evaluated via the closed-form approximation~(4) with the adjusted 
coefficients of variation $\gamma_X/\sqrt{n}$ and $\gamma_Y/\sqrt{n}$ 
(see~(9)). This approximation, originally proposed by 
Celano and Castagliola~(2016), is used as a fast surrogate for the 
Markov-chain construction. The exact density~(2) of Nadarajah and 
Okorie~(2016) is retained solely for high-accuracy validation of the 
reported TARL values (see Subsection 5.6).

Let $\mathbf{Q}=(Q_{ij})\in[0,1]^{m\times m}$ and $\mathbf{p}^{\mathrm{abs}}=(p_i^{\mathrm{abs}})\in[0,1]^m$, and let $\pmb{\pi}_0$ be the initial distribution over the safe states, taken as a unit mass on the state whose centre is closest to $z_0$. Define $\pmb{\alpha}_k = \pmb{\pi}_0\,\mathbf{Q}^k$ as the (sub-)distribution over the safe states after $k$ transitions. Then the probability of a signal exactly at inspection $k$ is
\begin{align}\label{eq:pk}
	\Pr(T=k) = \pmb{\alpha}_{k-1}^{\,T}\,\mathbf{p}^{\mathrm{abs}},
	\qquad k=1,\dots,I,
\end{align}
and the probability of no signal during the run is $\Pr(T>I) = \mathbf{1}^T\,\pmb{\alpha}_I$. Substituting into \eqref{eq:TARL} yields
\begin{align}\label{eq:TARL_MC}
	\mathrm{TARL}(I)
	= \sum_{k=1}^{I} k\,\pmb{\alpha}_{k-1}^{\,T}\,\mathbf{p}^{\mathrm{abs}}
	+ (I+1)\,\mathbf{1}^T\,\pmb{\alpha}_I.
\end{align}

\subsection{UCL calibration}

The design problem \eqref{eq:UCL_problem} reduces to a one-dimensional root-finding problem in $\mathit{UCL}$. Define
\begin{align*}
	g(\mathit{UCL}) = \mathrm{TARL}_0(\mathit{UCL}) - \mathit{TARL}_0^*.
\end{align*}
Since $\mathrm{TARL}_0$ is monotone increasing in $\mathit{UCL}$ (a wider safe region yields a longer in-control run length), and bounded by $I+1$, the equation $g(\mathit{UCL})=0$ admits a unique solution whenever the target is feasible. We use Brent's method on the bracket $[\mathit{UCL}_{\min},\mathit{UCL}_{\max}]$, where $\mathit{UCL}_{\min}=\max\{z_0,\,F^{-1}_{\hat Z}(0.5)\}$ and $\mathit{UCL}_{\max}=F^{-1}_{\hat Z}(0.999)$, with automatic bracket widening when the initial guess fails to enclose the root.

\subsection{Implementation considerations}

A naive implementation of the Markov-chain build evaluates $F_{\hat Z}$ at every pair $(i,j)$, requiring $\mathcal{O}(m^2)$ function calls per UCL evaluation. By stacking the $m+1$ grid edges into a row vector and the $m$ state centres into a column vector, the matrix of $z$-edges $z_{i,j}=(c_j-h/2-(1-\lambda)c_i)/\lambda$ can be assembled in a single \texttt{numpy} broadcast, and $F_{\hat Z}$ can then be applied element-wise in $\mathcal{O}(m^2)$ \emph{vectorized} operations. This optimization, together with the use of \texttt{scipy.stats.norm} for the standard-normal c.d.f.\ and quantile function, reduces the time to build the entire $(\lambda,\gamma_X,\gamma_Y,\rho_0,n)$ factorial table from hours to a few seconds. Numerical leakage caused by tail mass falling below $w_{\min}$ is reabsorbed into the first state to keep each row of $\mathbf{Q}$ a probability distribution.

We use $m=80$ states in the numerical results of Section~\ref{sec:numerical}, which is consistent with the value used by \citet{Tran2021} and is sufficient for an accuracy of $10^{-3}$ on the TARL values. Larger values of $m$ (e.g., $m=150$ or $m=200$) yield essentially identical numbers at a modest computational cost.

%%%%%%%%%%%%%%%%%%%%%%%%%%%%%%%%%%%%%%%%%%%%%%%%%%%%%%%%%%%%%%%%%%%%%%%%%%%%%%%%
\section{Numerical Analysis}
\label{sec:numerical}

This section assesses the statistical performance of the proposed EWMA-RZ$^+$ chart through a comprehensive numerical study based on the Markov-chain construction of Section~\ref{sec:optimization}. The objectives are threefold: (i) to verify that the calibration procedure produces upper control limits achieving the prescribed in-control $TARL_0$ across a broad range of process configurations; (ii) to characterize the out-of-control performance $TARL_1$ under upward and downward shifts of the in-control ratio $z_0$; and (iii) to identify how the design parameters $(\lambda,\gamma_X,\gamma_Y,\rho_0,n,I)$ jointly affect detection speed. To this end, we conduct a full factorial experiment whose grid is described in Section~\ref{ssec:factorial}, report the calibrated control limits in Section~\ref{ssec:limits}, and analyze the detection profiles in Section~\ref{ssec:detection}. A direct comparison with the Shewhart-type ShRZ chart of \citet{Tran2021} is provided in Section~\ref{ssec:comparison}. Section~\ref{ssec:robustness} closes the numerical analysis with an empirical robustness study against departures from bivariate normality. All numerical computations were performed with $m = 80$ Markov-chain states and Brent's root-finder tolerance $10^{-5}$, which we verified are sufficient to ensure four-decimal stability of the reported $TARL$ values.

\subsection{Factorial design of the study}
\label{ssec:factorial}

The behaviour of the proposed EWMA-RZ$^+$ chart is investigated through an extensive factorial study, with the parameter grid summarized in \tab{tab:design}. The grid matches the one used in \citet{Tran2021} for the ShRZ chart, so that the EWMA-RZ$^+$ and ShRZ charts can be compared one-to-one. The in-control value of the ratio is normalized to $z_0=1$ without loss of generality, since changing $z_0$ only rescales the chart through the in-control standard-deviation ratio $\omega_0 = z_0\,\gamma_X/\gamma_Y$.

\begin{table}[!htbp]
	\centering
	\caption{Factorial design used in the numerical analysis.}\label{tab:design}
	\begin{tabular}{l l}
		\toprule
		Parameter & Levels \\
		\midrule
		Horizon $I$                      & $10,\;30$ \\
		Target $\mathit{TARL}_0^*$       & $10,\;30$ (chosen equal to $I$) \\
		Smoothing constant $\lambda$     & $0.1,\;0.2$ \\
		Sample size $n$                  & $1,\;5,\;7,\;10,\;15$ \\
		In-control correlation $\rho_0$  & $-0.8,\;-0.4,\;0.0,\;0.4,\;0.8$ \\
		CV pair (equal)                  & $(0.01,0.01),\;(0.2,0.2)$ \\
		CV pair (unequal)                & $(0.01,0.2),\;(0.2,0.01)$ \\
		Shift factor $\tau$              & $0.90,\,0.95,\,0.98,\,0.99,\,1.00,\,1.01,\,1.02,\,1.05,\,1.10$ \\
		Markov-chain states $m$          & $80$ \\
		\bottomrule
	\end{tabular}
\end{table}

For each combination of the in-control parameters $(\lambda,\gamma_X,\gamma_Y,\rho_0,n,I,\mathit{TARL}_0^*)$, the upper control limit $\mathit{UCL}$ is calibrated by solving \eqref{eq:UCL_problem}, then $\mathrm{TARL}_1$ is evaluated at the same $\mathit{UCL}$ across the nine values of $\tau$.

\subsection{Organization of the results}

The complete set of numerical results is reported in two groups of tables. Group A (Tables~\ref{tab:limits_eq_T10}--\ref{tab:limits_uneq_T30}) contains the \emph{design tables}: for each combination of $(\lambda,\gamma_X,\gamma_Y,\rho_0,n,I)$, we report the calibrated upper control limit $\mathit{UCL}$ and the achieved in-control $\mathrm{TARL}_0$. Group B (Tables~\ref{tab:tarl1_equal_T10_lam01}--\ref{tab:tarl1_unequal_T30_lam02}) contains the \emph{detection tables}: at each calibrated $\mathit{UCL}$, the out-of-control $\mathrm{TARL}_1$ is evaluated across the nine values of $\tau$. To keep each detection table within page width, Group B is split by $\lambda\in\{0.1,0.2\}$. In total, we report $4$ design tables and $8$ detection tables, covering all $2\times 2\times 5\times 5\times 2 = 200$ design configurations and the $200\times 9 = 1{,}800$ detection scenarios of the factorial study described in \tab{tab:design}.

\subsection{Calibrated control limits and in-control TARL (Group A)}
\label{ssec:limits}

Tables~\ref{tab:limits_eq_T10}--\ref{tab:limits_uneq_T30} report the calibrated $\mathit{UCL}$ and the achieved $\mathrm{TARL}_0$ for every combination of the in-control parameters. Three patterns emerge consistently across the four tables.

\paragraph{Effect of the sample size $n$.} For fixed $(\lambda,\gamma_X,\gamma_Y,\rho_0)$, the calibrated $\mathit{UCL}$ decreases monotonically with $n$. The reason is that a larger sample reduces the variance of $\hat Z_i$ (the sample-mean CVs scale as $1/\sqrt{n}$, see eq.~\eqref{eq:Fzhat}) and therefore tightens the in-control distribution of the EWMA statistic. A tighter $\mathit{UCL}$ is then needed to preserve the same in-control TARL. For example, in the equal-CV case $(\gamma_X,\gamma_Y)=(0.2,0.2)$ with $\lambda=0.2$, $\rho_0=0.0$, $I=10$ (Table~\ref{tab:limits_eq_T10}), $\mathit{UCL}$ falls from $1.2134$ at $n=1$ to $1.0446$ at $n=15$.

\paragraph{Effect of the in-control correlation $\rho_0$.} For fixed $(\lambda,\gamma_X,\gamma_Y,n)$, the calibrated $\mathit{UCL}$ is monotonically decreasing in $\rho_0$. A positive correlation between $X$ and $Y$ causes them to move together, which dampens the variance of the ratio $\hat Z_i$; the chart can therefore afford a tighter limit. In the same equal-CV case as above, $\mathit{UCL}$ at $n=5$ drops from $1.1078$ at $\rho_0=-0.8$ to $1.0330$ at $\rho_0=+0.8$.

\paragraph{Effect of the smoothing constant $\lambda$.} For fixed $(\gamma_X,\gamma_Y,\rho_0,n)$, $\lambda=0.2$ yields a larger $\mathit{UCL}$ than $\lambda=0.1$, by an approximately constant factor of $\sqrt{(2-0.1)/(2-0.2)\cdot 0.2/0.1}\approx \sqrt{2.11}\approx 1.45$. This is the well-known steady-state scaling of the EWMA recursion: the standard deviation of the EWMA grows with $\sqrt{\lambda/(2-\lambda)}$ at the long-run limit, although for the very short runs considered here the transient term plays a non-negligible role as well.

\paragraph{Accuracy of the calibration.} Across all four design tables, the achieved $\mathrm{TARL}_0$ matches the target $\mathit{TARL}_0^*$ to within a few thousandths of a unit. This confirms that the Markov-chain build with $m=80$ states and the Brent root-finder with tolerance $10^{-5}$ are sufficient for the design problem; tighter tolerances and larger $m$ produce essentially the same numbers.

% ====== include the full LIMITS tables ======
\begin{longtable}{c c r r r r}
\caption{Calibrated upper control limits $UCL$ and achieved in-control $TARL_0$ for the EWMA-RZ$^+$ chart under equal coefficients of variation ($\gamma_X = \gamma_Y$), horizon $I = 10$, target $TARL_0 = 10$, $z_0 = 1$.}\label{tab:limits_eq_T10}\\
\toprule
$\rho_0$ & $n$ & \multicolumn{2}{c}{$\lambda = 0.1$} & \multicolumn{2}{c}{$\lambda = 0.2$} \\
\cmidrule(lr){3-4}\cmidrule(lr){5-6}
 & & $UCL$ & $TARL_0$ & $UCL$ & $TARL_0$ \\
\midrule
\endfirsthead
\multicolumn{6}{l}{\textit{Table \ref{tab:limits_eq_T10} (continued)}}\\
\toprule
$\rho_0$ & $n$ & \multicolumn{2}{c}{$\lambda = 0.1$} & \multicolumn{2}{c}{$\lambda = 0.2$} \\
\cmidrule(lr){3-4}\cmidrule(lr){5-6}
 & & $UCL$ & $TARL_0$ & $UCL$ & $TARL_0$ \\
\midrule
\endhead
\bottomrule
\endlastfoot
\multicolumn{6}{l}{\textit{$\gamma_X = 0.01$, $\gamma_Y = 0.01$}} \\
-0.8 & 1 & 1.0052 & 10.000 & 1.0101 & 10.000 \\
 & 5 & 1.0022 & 10.001 & 1.0044 & 10.000 \\
 & 7 & 1.0018 & 10.001 & 1.0037 & 10.000 \\
 & 10 & 1.0015 & 10.001 & 1.0031 & 10.000 \\
 & 15 & 1.0012 & 10.001 & 1.0025 & 9.999 \\
\addlinespace[0.2em]
-0.4 & 1 & 1.0046 & 10.000 & 1.0089 & 10.000 \\
 & 5 & 1.0019 & 10.001 & 1.0039 & 10.000 \\
 & 7 & 1.0016 & 10.001 & 1.0033 & 10.000 \\
 & 10 & 1.0013 & 10.001 & 1.0027 & 10.000 \\
 & 15 & 1.0011 & 10.001 & 1.0022 & 9.999 \\
\addlinespace[0.2em]
0.0 & 1 & 1.0038 & 10.000 & 1.0074 & 10.000 \\
 & 5 & 1.0016 & 10.001 & 1.0033 & 10.000 \\
 & 7 & 1.0014 & 10.001 & 1.0028 & 10.000 \\
 & 10 & 1.0011 & 10.001 & 1.0023 & 9.999 \\
 & 15 & 1.0009 & 10.001 & 1.0019 & 9.998 \\
\addlinespace[0.2em]
0.4 & 1 & 1.0029 & 10.001 & 1.0057 & 10.000 \\
 & 5 & 1.0012 & 10.001 & 1.0025 & 9.999 \\
 & 7 & 1.0010 & 10.001 & 1.0021 & 9.999 \\
 & 10 & 1.0009 & 10.001 & 1.0018 & 9.998 \\
 & 15 & 1.0007 & 10.001 & 1.0014 & 10.003 \\
\addlinespace[0.2em]
0.8 & 1 & 1.0016 & 10.001 & 1.0033 & 10.000 \\
 & 5 & 1.0007 & 10.001 & 1.0014 & 10.003 \\
 & 7 & 1.0006 & 10.001 & 1.0012 & 10.003 \\
 & 10 & 1.0005 & 10.001 & 1.0010 & 10.003 \\
 & 15 & 1.0004 & 10.001 & 1.0008 & 10.003 \\
\addlinespace[0.2em]
\midrule
\multicolumn{6}{l}{\textit{$\gamma_X = 0.2$, $\gamma_Y = 0.2$}} \\
-0.8 & 1 & 1.1643 & 10.000 & 1.3305 & 10.000 \\
 & 5 & 1.0618 & 10.327 & 1.1078 & 10.000 \\
 & 7 & 1.0541 & 9.678 & 1.0849 & 10.000 \\
 & 10 & 1.0469 & 9.855 & 1.0712 & 10.121 \\
 & 15 & 1.0392 & 10.000 & 1.0601 & 9.875 \\
\addlinespace[0.2em]
-0.4 & 1 & 1.1579 & 9.803 & 1.2702 & 10.079 \\
 & 5 & 1.0566 & 9.650 & 1.0901 & 10.000 \\
 & 7 & 1.0493 & 9.822 & 1.0749 & 10.086 \\
 & 10 & 1.0424 & 9.968 & 1.0645 & 9.818 \\
 & 15 & 1.0330 & 10.000 & 1.0542 & 9.983 \\
\addlinespace[0.2em]
0.0 & 1 & 1.1234 & 10.000 & 1.2134 & 10.000 \\
 & 5 & 1.0500 & 9.848 & 1.0759 & 10.111 \\
 & 7 & 1.0432 & 9.979 & 1.0657 & 9.831 \\
 & 10 & 1.0344 & 10.000 & 1.0562 & 9.974 \\
 & 15 & 1.0263 & 10.000 & 1.0446 & 10.000 \\
\addlinespace[0.2em]
0.4 & 1 & 1.0878 & 10.257 & 1.1603 & 10.000 \\
 & 5 & 1.0385 & 10.000 & 1.0621 & 9.951 \\
 & 7 & 1.0307 & 10.000 & 1.0518 & 10.000 \\
 & 10 & 1.0244 & 10.000 & 1.0419 & 10.000 \\
 & 15 & 1.0188 & 10.000 & 1.0332 & 10.000 \\
\addlinespace[0.2em]
0.8 & 1 & 1.0477 & 10.000 & 1.0859 & 10.000 \\
 & 5 & 1.0182 & 10.000 & 1.0330 & 10.000 \\
 & 7 & 1.0149 & 10.000 & 1.0273 & 10.000 \\
 & 10 & 1.0121 & 10.000 & 1.0224 & 10.000 \\
 & 15 & 1.0096 & 10.000 & 1.0180 & 10.000 \\
\addlinespace[0.2em]
\end{longtable}

\begin{longtable}{c c r r r r}
\caption{Calibrated upper control limits $UCL$ and achieved in-control $TARL_0$ for the EWMA-RZ$^+$ chart under unequal coefficients of variation ($\gamma_X \neq \gamma_Y$), horizon $I = 10$, target $TARL_0 = 10$, $z_0 = 1$.}\label{tab:limits_uneq_T10}\\
\toprule
$\rho_0$ & $n$ & \multicolumn{2}{c}{$\lambda = 0.1$} & \multicolumn{2}{c}{$\lambda = 0.2$} \\
\cmidrule(lr){3-4}\cmidrule(lr){5-6}
 & & $UCL$ & $TARL_0$ & $UCL$ & $TARL_0$ \\
\midrule
\endfirsthead
\multicolumn{6}{l}{\textit{Table \ref{tab:limits_uneq_T10} (continued)}}\\
\toprule
$\rho_0$ & $n$ & \multicolumn{2}{c}{$\lambda = 0.1$} & \multicolumn{2}{c}{$\lambda = 0.2$} \\
\cmidrule(lr){3-4}\cmidrule(lr){5-6}
 & & $UCL$ & $TARL_0$ & $UCL$ & $TARL_0$ \\
\midrule
\endhead
\bottomrule
\endlastfoot
\multicolumn{6}{l}{\textit{$\gamma_X = 0.01$, $\gamma_Y = 0.2$}} \\
-0.8 & 1 & 1.0993 & 10.000 & 1.1800 & 10.073 \\
 & 5 & 1.0333 & 10.273 & 1.0605 & 10.000 \\
 & 7 & 1.0293 & 9.611 & 1.0475 & 10.000 \\
 & 10 & 1.0254 & 9.803 & 1.0386 & 10.068 \\
 & 15 & 1.0215 & 9.973 & 1.0327 & 9.824 \\
\addlinespace[0.2em]
-0.4 & 1 & 1.0973 & 10.000 & 1.1767 & 10.074 \\
 & 5 & 1.0327 & 10.273 & 1.0594 & 10.000 \\
 & 7 & 1.0287 & 9.612 & 1.0466 & 10.000 \\
 & 10 & 1.0249 & 9.803 & 1.0379 & 10.068 \\
 & 15 & 1.0211 & 9.973 & 1.0320 & 9.824 \\
\addlinespace[0.2em]
0.0 & 1 & 1.0954 & 10.000 & 1.1733 & 10.075 \\
 & 5 & 1.0321 & 10.274 & 1.0582 & 10.000 \\
 & 7 & 1.0282 & 9.612 & 1.0457 & 10.000 \\
 & 10 & 1.0244 & 9.803 & 1.0372 & 10.069 \\
 & 15 & 1.0207 & 9.973 & 1.0314 & 9.825 \\
\addlinespace[0.2em]
0.4 & 1 & 1.0935 & 10.000 & 1.1698 & 10.075 \\
 & 5 & 1.0314 & 10.274 & 1.0571 & 10.000 \\
 & 7 & 1.0276 & 9.612 & 1.0448 & 10.000 \\
 & 10 & 1.0240 & 9.803 & 1.0364 & 10.069 \\
 & 15 & 1.0203 & 9.973 & 1.0308 & 9.825 \\
\addlinespace[0.2em]
0.8 & 1 & 1.0916 & 10.000 & 1.1661 & 10.073 \\
 & 5 & 1.0307 & 10.273 & 1.0559 & 10.000 \\
 & 7 & 1.0270 & 9.611 & 1.0439 & 10.000 \\
 & 10 & 1.0234 & 9.803 & 1.0356 & 10.068 \\
 & 15 & 1.0198 & 9.973 & 1.0301 & 9.824 \\
\addlinespace[0.2em]
\midrule
\multicolumn{6}{l}{\textit{$\gamma_X = 0.2$, $\gamma_Y = 0.01$}} \\
-0.8 & 1 & 1.0565 & 10.000 & 1.1102 & 10.000 \\
 & 5 & 1.0239 & 10.000 & 1.0481 & 10.000 \\
 & 7 & 1.0200 & 10.000 & 1.0406 & 10.000 \\
 & 10 & 1.0166 & 10.000 & 1.0339 & 10.000 \\
 & 15 & 1.0135 & 10.000 & 1.0276 & 10.000 \\
\addlinespace[0.2em]
-0.4 & 1 & 1.0529 & 10.000 & 1.1061 & 10.000 \\
 & 5 & 1.0229 & 10.000 & 1.0468 & 10.000 \\
 & 7 & 1.0193 & 10.000 & 1.0395 & 10.000 \\
 & 10 & 1.0161 & 10.000 & 1.0330 & 10.000 \\
 & 15 & 1.0131 & 10.000 & 1.0269 & 10.000 \\
\addlinespace[0.2em]
0.0 & 1 & 1.0491 & 10.000 & 1.1019 & 10.000 \\
 & 5 & 1.0219 & 10.000 & 1.0455 & 10.000 \\
 & 7 & 1.0185 & 10.000 & 1.0384 & 10.000 \\
 & 10 & 1.0155 & 10.000 & 1.0322 & 10.000 \\
 & 15 & 1.0126 & 10.000 & 1.0263 & 10.000 \\
\addlinespace[0.2em]
0.4 & 1 & 1.0468 & 10.027 & 1.0977 & 10.000 \\
 & 5 & 1.0209 & 10.000 & 1.0441 & 10.000 \\
 & 7 & 1.0177 & 10.000 & 1.0374 & 10.000 \\
 & 10 & 1.0149 & 10.000 & 1.0313 & 10.000 \\
 & 15 & 1.0122 & 10.000 & 1.0256 & 10.000 \\
\addlinespace[0.2em]
0.8 & 1 & 1.0466 & 10.090 & 1.0944 & 10.017 \\
 & 5 & 1.0205 & 10.024 & 1.0428 & 10.000 \\
 & 7 & 1.0173 & 10.016 & 1.0363 & 10.000 \\
 & 10 & 1.0144 & 10.009 & 1.0304 & 10.000 \\
 & 15 & 1.0118 & 10.002 & 1.0249 & 10.000 \\
\addlinespace[0.2em]
\end{longtable}

\begin{longtable}{c c r r r r}
\caption{Calibrated upper control limits $UCL$ and achieved in-control $TARL_0$ for the EWMA-RZ$^+$ chart under equal coefficients of variation ($\gamma_X = \gamma_Y$), horizon $I = 30$, target $TARL_0 = 30$, $z_0 = 1$.}\label{tab:limits_eq_T30}\\
\toprule
$\rho_0$ & $n$ & \multicolumn{2}{c}{$\lambda = 0.1$} & \multicolumn{2}{c}{$\lambda = 0.2$} \\
\cmidrule(lr){3-4}\cmidrule(lr){5-6}
 & & $UCL$ & $TARL_0$ & $UCL$ & $TARL_0$ \\
\midrule
\endfirsthead
\multicolumn{6}{l}{\textit{Table \ref{tab:limits_eq_T30} (continued)}}\\
\toprule
$\rho_0$ & $n$ & \multicolumn{2}{c}{$\lambda = 0.1$} & \multicolumn{2}{c}{$\lambda = 0.2$} \\
\cmidrule(lr){3-4}\cmidrule(lr){5-6}
 & & $UCL$ & $TARL_0$ & $UCL$ & $TARL_0$ \\
\midrule
\endhead
\bottomrule
\endlastfoot
\multicolumn{6}{l}{\textit{$\gamma_X = 0.01$, $\gamma_Y = 0.01$}} \\
-0.8 & 1 & 1.0116 & 30.000 & 1.0175 & 30.000 \\
 & 5 & 1.0051 & 29.999 & 1.0077 & 29.999 \\
 & 7 & 1.0043 & 29.999 & 1.0065 & 29.998 \\
 & 10 & 1.0036 & 29.999 & 1.0054 & 29.997 \\
 & 15 & 1.0029 & 29.999 & 1.0044 & 30.003 \\
\addlinespace[0.2em]
-0.4 & 1 & 1.0102 & 30.000 & 1.0154 & 30.000 \\
 & 5 & 1.0045 & 29.999 & 1.0068 & 29.998 \\
 & 7 & 1.0038 & 29.999 & 1.0057 & 29.998 \\
 & 10 & 1.0032 & 29.999 & 1.0048 & 30.003 \\
 & 15 & 1.0026 & 29.999 & 1.0039 & 30.003 \\
\addlinespace[0.2em]
0.0 & 1 & 1.0086 & 29.999 & 1.0129 & 30.000 \\
 & 5 & 1.0038 & 29.999 & 1.0057 & 29.998 \\
 & 7 & 1.0032 & 29.999 & 1.0048 & 30.003 \\
 & 10 & 1.0027 & 29.999 & 1.0040 & 30.003 \\
 & 15 & 1.0022 & 29.999 & 1.0033 & 30.003 \\
\addlinespace[0.2em]
0.4 & 1 & 1.0066 & 29.999 & 1.0100 & 30.000 \\
 & 5 & 1.0029 & 29.999 & 1.0044 & 30.003 \\
 & 7 & 1.0025 & 29.999 & 1.0037 & 30.003 \\
 & 10 & 1.0021 & 29.999 & 1.0031 & 30.003 \\
 & 15 & 1.0017 & 29.999 & 1.0026 & 30.003 \\
\addlinespace[0.2em]
0.8 & 1 & 1.0038 & 29.999 & 1.0057 & 29.998 \\
 & 5 & 1.0017 & 29.999 & 1.0026 & 30.003 \\
 & 7 & 1.0014 & 29.999 & 1.0022 & 30.006 \\
 & 10 & 1.0012 & 29.999 & 1.0018 & 30.006 \\
 & 15 & 1.0010 & 29.999 & 1.0015 & 30.006 \\
\addlinespace[0.2em]
\midrule
\multicolumn{6}{l}{\textit{$\gamma_X = 0.2$, $\gamma_Y = 0.2$}} \\
-0.8 & 1 & 1.3647 & 30.000 & 1.6058 & 30.000 \\
 & 5 & 1.1209 & 30.000 & 1.1836 & 30.000 \\
 & 7 & 1.0973 & 30.000 & 1.1483 & 30.000 \\
 & 10 & 1.0780 & 30.000 & 1.1219 & 29.934 \\
 & 15 & 1.0613 & 30.000 & 1.0972 & 30.000 \\
\addlinespace[0.2em]
-0.4 & 1 & 1.3107 & 30.000 & 1.5117 & 30.000 \\
 & 5 & 1.1030 & 30.000 & 1.1571 & 30.000 \\
 & 7 & 1.0834 & 30.000 & 1.1282 & 30.023 \\
 & 10 & 1.0672 & 30.000 & 1.1065 & 30.000 \\
 & 15 & 1.0542 & 30.092 & 1.0843 & 30.000 \\
\addlinespace[0.2em]
0.0 & 1 & 1.2531 & 30.000 & 1.4016 & 30.000 \\
 & 5 & 1.0837 & 30.000 & 1.1299 & 30.066 \\
 & 7 & 1.0682 & 30.000 & 1.1082 & 30.000 \\
 & 10 & 1.0562 & 30.075 & 1.0879 & 30.000 \\
 & 15 & 1.0470 & 29.754 & 1.0700 & 30.000 \\
\addlinespace[0.2em]
0.4 & 1 & 1.1802 & 30.180 & 1.2828 & 30.000 \\
 & 5 & 1.0621 & 30.029 & 1.0984 & 30.000 \\
 & 7 & 1.0533 & 30.198 & 1.0809 & 30.000 \\
 & 10 & 1.0452 & 29.854 & 1.0662 & 30.000 \\
 & 15 & 1.0373 & 30.000 & 1.0530 & 30.000 \\
\addlinespace[0.2em]
0.8 & 1 & 1.0971 & 29.802 & 1.1461 & 30.000 \\
 & 5 & 1.0374 & 30.000 & 1.0535 & 30.000 \\
 & 7 & 1.0311 & 30.000 & 1.0445 & 30.000 \\
 & 10 & 1.0256 & 30.000 & 1.0373 & 30.090 \\
 & 15 & 1.0206 & 30.000 & 1.0306 & 29.966 \\
\addlinespace[0.2em]
\end{longtable}

\begin{longtable}{c c r r r r}
\caption{Calibrated upper control limits $UCL$ and achieved in-control $TARL_0$ for the EWMA-RZ$^+$ chart under unequal coefficients of variation ($\gamma_X \neq \gamma_Y$), horizon $I = 30$, target $TARL_0 = 30$, $z_0 = 1$.}\label{tab:limits_uneq_T30}\\
\toprule
$\rho_0$ & $n$ & \multicolumn{2}{c}{$\lambda = 0.1$} & \multicolumn{2}{c}{$\lambda = 0.2$} \\
\cmidrule(lr){3-4}\cmidrule(lr){5-6}
 & & $UCL$ & $TARL_0$ & $UCL$ & $TARL_0$ \\
\midrule
\endfirsthead
\multicolumn{6}{l}{\textit{Table \ref{tab:limits_uneq_T30} (continued)}}\\
\toprule
$\rho_0$ & $n$ & \multicolumn{2}{c}{$\lambda = 0.1$} & \multicolumn{2}{c}{$\lambda = 0.2$} \\
\cmidrule(lr){3-4}\cmidrule(lr){5-6}
 & & $UCL$ & $TARL_0$ & $UCL$ & $TARL_0$ \\
\midrule
\endhead
\bottomrule
\endlastfoot
\multicolumn{6}{l}{\textit{$\gamma_X = 0.01$, $\gamma_Y = 0.2$}} \\
-0.8 & 1 & 1.2089 & 30.000 & 1.3444 & 30.000 \\
 & 5 & 1.0673 & 30.000 & 1.1020 & 30.000 \\
 & 7 & 1.0540 & 30.000 & 1.0823 & 30.000 \\
 & 10 & 1.0433 & 30.000 & 1.0661 & 30.000 \\
 & 15 & 1.0339 & 30.000 & 1.0537 & 30.000 \\
\addlinespace[0.2em]
-0.4 & 1 & 1.2049 & 30.000 & 1.3378 & 30.000 \\
 & 5 & 1.0660 & 30.000 & 1.1001 & 30.000 \\
 & 7 & 1.0530 & 30.000 & 1.0807 & 30.000 \\
 & 10 & 1.0424 & 30.000 & 1.0649 & 30.000 \\
 & 15 & 1.0333 & 30.000 & 1.0527 & 30.000 \\
\addlinespace[0.2em]
0.0 & 1 & 1.2009 & 30.000 & 1.3312 & 30.000 \\
 & 5 & 1.0648 & 30.000 & 1.0982 & 30.000 \\
 & 7 & 1.0520 & 30.000 & 1.0791 & 30.000 \\
 & 10 & 1.0416 & 30.000 & 1.0636 & 30.000 \\
 & 15 & 1.0326 & 30.000 & 1.0517 & 30.000 \\
\addlinespace[0.2em]
0.4 & 1 & 1.1969 & 30.000 & 1.3246 & 30.000 \\
 & 5 & 1.0635 & 30.000 & 1.0962 & 30.000 \\
 & 7 & 1.0510 & 30.000 & 1.0776 & 30.000 \\
 & 10 & 1.0408 & 30.000 & 1.0624 & 30.000 \\
 & 15 & 1.0320 & 30.000 & 1.0506 & 30.000 \\
\addlinespace[0.2em]
0.8 & 1 & 1.1929 & 30.000 & 1.3179 & 30.000 \\
 & 5 & 1.0621 & 30.000 & 1.0942 & 30.000 \\
 & 7 & 1.0499 & 30.000 & 1.0759 & 30.000 \\
 & 10 & 1.0399 & 30.000 & 1.0611 & 30.000 \\
 & 15 & 1.0313 & 30.000 & 1.0496 & 30.000 \\
\addlinespace[0.2em]
\midrule
\multicolumn{6}{l}{\textit{$\gamma_X = 0.2$, $\gamma_Y = 0.01$}} \\
-0.8 & 1 & 1.1270 & 30.000 & 1.1910 & 30.000 \\
 & 5 & 1.0557 & 30.000 & 1.0844 & 30.000 \\
 & 7 & 1.0470 & 30.000 & 1.0712 & 30.000 \\
 & 10 & 1.0392 & 30.000 & 1.0595 & 30.000 \\
 & 15 & 1.0320 & 30.000 & 1.0485 & 30.000 \\
\addlinespace[0.2em]
-0.4 & 1 & 1.1227 & 30.000 & 1.1856 & 30.000 \\
 & 5 & 1.0543 & 30.000 & 1.0825 & 30.000 \\
 & 7 & 1.0458 & 30.000 & 1.0696 & 30.000 \\
 & 10 & 1.0383 & 30.000 & 1.0582 & 30.000 \\
 & 15 & 1.0312 & 30.000 & 1.0475 & 30.000 \\
\addlinespace[0.2em]
0.0 & 1 & 1.1183 & 30.000 & 1.1802 & 30.000 \\
 & 5 & 1.0529 & 30.000 & 1.0805 & 30.000 \\
 & 7 & 1.0447 & 30.000 & 1.0680 & 30.000 \\
 & 10 & 1.0374 & 30.000 & 1.0569 & 30.000 \\
 & 15 & 1.0305 & 30.000 & 1.0465 & 30.000 \\
\addlinespace[0.2em]
0.4 & 1 & 1.1138 & 30.000 & 1.1747 & 30.000 \\
 & 5 & 1.0514 & 30.000 & 1.0785 & 30.000 \\
 & 7 & 1.0435 & 30.000 & 1.0664 & 30.000 \\
 & 10 & 1.0364 & 30.000 & 1.0556 & 30.000 \\
 & 15 & 1.0298 & 30.000 & 1.0454 & 30.000 \\
\addlinespace[0.2em]
0.8 & 1 & 1.1092 & 30.000 & 1.1692 & 30.000 \\
 & 5 & 1.0499 & 30.000 & 1.0765 & 30.000 \\
 & 7 & 1.0422 & 30.000 & 1.0647 & 30.000 \\
 & 10 & 1.0354 & 30.000 & 1.0542 & 30.000 \\
 & 15 & 1.0290 & 30.000 & 1.0443 & 30.000 \\
\addlinespace[0.2em]
\end{longtable}

\subsection{Detection performance (Group B)}
\label{ssec:detection}

\begin{figure}[H]
	\centering
	\includegraphics[width=0.85\linewidth]{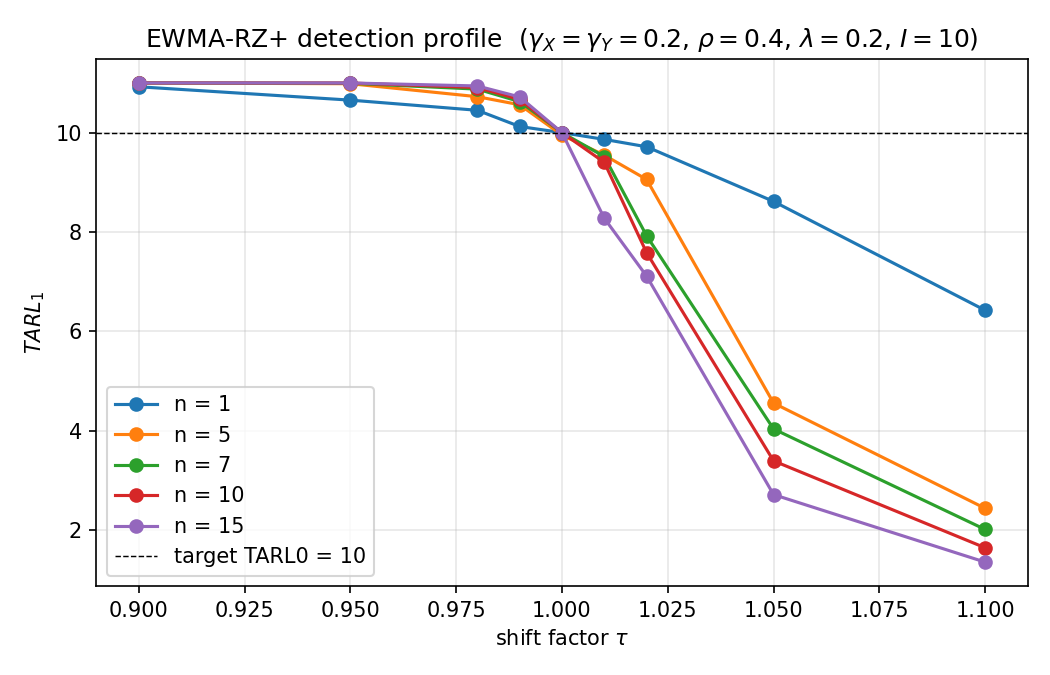}
	\caption{Detection profile of the EWMA-RZ$^+$ chart: $\mathrm{TARL}_1$ as a function of the shift factor $\tau$, for the equal-CV case $(\gamma_X,\gamma_Y)=(0.2,0.2)$, $\rho_0=0.4$, $\lambda=0.2$, $I=10$, $\mathit{TARL}_0^*=10$.}
	\label{fig:profile}
\end{figure}

Tables~\ref{tab:tarl1_equal_T10_lam01}--\ref{tab:tarl1_unequal_T30_lam02} report the out-of-control $\mathrm{TARL}_1$ values across the nine shift levels $\tau$. We summarize the most salient findings below; the figure-level visualization of one slice is shown in \fig{fig:profile}.

\paragraph{Asymmetry in $\tau$.} As expected for an upper-sided chart, downward shifts ($\tau<1$) are essentially undetected: $\mathrm{TARL}_1$ values are close to $I+1$ across every entry of the $\tau\in\{0.90,0.95\}$ columns. The chart's responsibility is to detect upward shifts ($\tau>1$); the downward-shift columns serve as a sanity check that the chart does not over-react in the wrong direction.

\paragraph{Effect of the sample size on detection.} For every upward shift, increasing $n$ reduces $\mathrm{TARL}_1$ monotonically. The effect is most pronounced for moderate shifts: in Table~\ref{tab:tarl1_equal_T10_lam02} (equal CV, $\lambda=0.2$, $I=10$), with $(\gamma_X,\gamma_Y)=(0.2,0.2)$ and $\rho_0=0.4$, $\mathrm{TARL}_1$ at $\tau=1.05$ drops from $8.615$ at $n=1$ to $2.713$ at $n=15$. At larger shifts ($\tau=1.10$), $\mathrm{TARL}_1$ approaches $1$ even for moderate $n$, which means the chart signals essentially at the first post-shift inspection.

\paragraph{Effect of the CV.} Small CVs sharpen the chart: when $(\gamma_X,\gamma_Y)=(0.01,0.01)$, even shifts as small as $\tau=1.01$ are detected within $1$--$2$ inspections regardless of $n$ (the $\mathit{UCL}$ is then very close to $z_0$). Large CVs disperse the chart: with $(\gamma_X,\gamma_Y)=(0.2,0.2)$, shifts below $\tau=1.05$ rarely yield $\mathrm{TARL}_1$ much below the target, especially for $n=1$. The unequal-CV cases $(0.01,0.2)$ and $(0.2,0.01)$ (Tables~\ref{tab:tarl1_unequal_T10_lam01}--\ref{tab:tarl1_unequal_T30_lam02}) are not symmetric: the chart is more sensitive when $\gamma_Y$ is the larger of the two, because $Y$ enters the denominator of $\hat Z_i$ and its variance is therefore amplified non-linearly through the ratio.

\paragraph{Effect of the correlation.} Once $\mathit{UCL}$ has been re-calibrated for the new $\rho_0$, the detection performance varies only mildly with the in-control correlation. The mild variation that remains is in favour of positive correlation: positive $\rho_0$ further reduces the variance of $\hat Z_i$ under the shift, which makes the EWMA statistic climb more deterministically toward $\mathit{UCL}$. The effect is at most a 5--10\% improvement in $\mathrm{TARL}_1$ when moving from $\rho_0=-0.8$ to $\rho_0=+0.8$ at fixed $n$ and $\tau$.

\paragraph{Effect of the smoothing constant.} Comparing the $\lambda=0.1$ and $\lambda=0.2$ versions of each detection table (e.g.\ Table~\ref{tab:tarl1_equal_T10_lam01} vs.\ Table~\ref{tab:tarl1_equal_T10_lam02}), we observe that $\lambda=0.1$ outperforms $\lambda=0.2$ for the smallest upward shifts ($\tau\in\{1.01,1.02\}$) by a few percent, while $\lambda=0.2$ is faster than $\lambda=0.1$ for $\tau\ge 1.05$ by a similarly small margin. For routine SPR monitoring, $\lambda=0.2$ is a robust default; if the user has a strong prior that the shift will be very small and persistent, $\lambda=0.1$ may be preferable. This trade-off is consistent with the long-run EWMA literature \citep{Lucas1990}.

\paragraph{Effect of the horizon $I$.} Comparing the $I=10$ tables with the $I=30$ tables, two trends stand out. First, the calibrated $\mathit{UCL}$ is larger for $I=30$ than for $I=10$ (the chart must afford a longer in-control run). Second, at fixed shift and sample size, $\mathrm{TARL}_1$ is roughly proportional to the horizon for very small shifts (e.g.\ $\tau=1.01$) but is essentially insensitive to the horizon for moderate or large shifts (e.g.\ $\tau\ge 1.05$). The chart effectively ``saturates'' for large shifts: it signals within $1$--$2$ inspections regardless of $I$.

% ====== include the full TARL1 tables ======
\begingroup
\setlength{\tabcolsep}{3pt}
\begin{longtable}{c c r r r r r r r r r}
\caption{$TARL_1$ values for the EWMA-RZ$^+$ chart with smoothing constant $\lambda = 0.1$, under equal coefficients of variation ($\gamma_X = \gamma_Y$), horizon $I = 10$, target $TARL_0 = 10$, $z_0 = 1$, across nine shift levels $\tau$.}\label{tab:tarl1_equal_T10_lam01}\\
\toprule
$\rho_0$ & $n$ & $\tau=0.9$ & $\tau=0.95$ & $\tau=0.98$ & $\tau=0.99$ & $\tau=1$ & $\tau=1.01$ & $\tau=1.02$ & $\tau=1.05$ & $\tau=1.1$ \\
\midrule
\endfirsthead
\multicolumn{11}{l}{\textit{Table \ref{tab:tarl1_equal_T10_lam01} (continued)}}\\
\toprule
$\rho_0$ & $n$ & $\tau=0.9$ & $\tau=0.95$ & $\tau=0.98$ & $\tau=0.99$ & $\tau=1$ & $\tau=1.01$ & $\tau=1.02$ & $\tau=1.05$ & $\tau=1.1$ \\
\midrule
\endhead
\bottomrule
\endlastfoot
\multicolumn{11}{l}{\textit{$\gamma_X = 0.01$, $\gamma_Y = 0.01$}} \\
-0.8 & 1 & 11.000 & 11.000 & 11.000 & 10.997 & 10.000 & 4.628 & 2.481 & 1.092 & 1.000 \\
 & 5 & 11.000 & 11.000 & 11.000 & 11.000 & 10.001 & 2.171 & 1.144 & 1.000 & 1.000 \\
 & 7 & 11.000 & 11.000 & 11.000 & 11.000 & 10.001 & 1.802 & 1.058 & 1.000 & 1.000 \\
 & 10 & 11.000 & 11.000 & 11.000 & 11.000 & 10.001 & 1.499 & 1.015 & 1.000 & 1.000 \\
 & 15 & 11.000 & 11.000 & 11.000 & 11.000 & 10.001 & 1.255 & 1.002 & 1.000 & 1.000 \\
\addlinespace[0.2em]
-0.4 & 1 & 11.000 & 11.000 & 11.000 & 10.998 & 10.000 & 4.238 & 2.156 & 1.043 & 1.000 \\
 & 5 & 11.000 & 11.000 & 11.000 & 11.000 & 10.001 & 1.888 & 1.075 & 1.000 & 1.000 \\
 & 7 & 11.000 & 11.000 & 11.000 & 11.000 & 10.001 & 1.579 & 1.024 & 1.000 & 1.000 \\
 & 10 & 11.000 & 11.000 & 11.000 & 11.000 & 10.001 & 1.336 & 1.004 & 1.000 & 1.000 \\
 & 15 & 11.000 & 11.000 & 11.000 & 11.000 & 10.001 & 1.151 & 1.000 & 1.000 & 1.000 \\
\addlinespace[0.2em]
0.0 & 1 & 11.000 & 11.000 & 11.000 & 11.000 & 10.000 & 3.693 & 1.790 & 1.011 & 1.000 \\
 & 5 & 11.000 & 11.000 & 11.000 & 11.000 & 10.001 & 1.579 & 1.024 & 1.000 & 1.000 \\
 & 7 & 11.000 & 11.000 & 11.000 & 11.000 & 10.001 & 1.347 & 1.005 & 1.000 & 1.000 \\
 & 10 & 11.000 & 11.000 & 11.000 & 11.000 & 10.001 & 1.176 & 1.000 & 1.000 & 1.000 \\
 & 15 & 11.000 & 11.000 & 11.000 & 11.000 & 10.001 & 1.062 & 1.000 & 1.000 & 1.000 \\
\addlinespace[0.2em]
0.4 & 1 & 11.000 & 11.000 & 11.000 & 11.000 & 10.001 & 2.877 & 1.386 & 1.000 & 1.000 \\
 & 5 & 11.000 & 11.000 & 11.000 & 11.000 & 10.001 & 1.255 & 1.002 & 1.000 & 1.000 \\
 & 7 & 11.000 & 11.000 & 11.000 & 11.000 & 10.001 & 1.124 & 1.000 & 1.000 & 1.000 \\
 & 10 & 11.000 & 11.000 & 11.000 & 11.000 & 10.001 & 1.044 & 1.000 & 1.000 & 1.000 \\
 & 15 & 11.000 & 11.000 & 11.000 & 11.000 & 10.001 & 1.008 & 1.000 & 1.000 & 1.000 \\
\addlinespace[0.2em]
0.8 & 1 & 11.000 & 11.000 & 11.000 & 11.000 & 10.001 & 1.580 & 1.024 & 1.000 & 1.000 \\
 & 5 & 11.000 & 11.000 & 11.000 & 11.000 & 10.001 & 1.008 & 1.000 & 1.000 & 1.000 \\
 & 7 & 11.000 & 11.000 & 11.000 & 11.000 & 10.001 & 1.001 & 1.000 & 1.000 & 1.000 \\
 & 10 & 11.000 & 11.000 & 11.000 & 11.000 & 10.001 & 1.000 & 1.000 & 1.000 & 1.000 \\
 & 15 & 11.000 & 11.000 & 11.000 & 11.000 & 10.001 & 1.000 & 1.000 & 1.000 & 1.000 \\
\addlinespace[0.2em]
\midrule
\multicolumn{11}{l}{\textit{$\gamma_X = 0.2$, $\gamma_Y = 0.2$}} \\
-0.8 & 1 & 10.838 & 10.580 & 10.140 & 10.072 & 10.000 & 9.315 & 9.210 & 8.875 & 7.022 \\
 & 5 & 10.999 & 10.938 & 10.583 & 10.466 & 10.327 & 9.202 & 8.911 & 5.514 & 3.954 \\
 & 7 & 11.000 & 10.971 & 10.705 & 10.595 & 9.678 & 9.381 & 9.050 & 5.423 & 3.649 \\
 & 10 & 11.000 & 10.997 & 10.799 & 10.697 & 9.855 & 9.523 & 9.140 & 5.238 & 3.258 \\
 & 15 & 11.000 & 10.999 & 10.961 & 10.777 & 10.000 & 9.616 & 6.772 & 4.887 & 2.754 \\
\addlinespace[0.2em]
-0.4 & 1 & 10.919 & 10.751 & 10.418 & 10.360 & 9.803 & 9.705 & 9.602 & 8.157 & 7.375 \\
 & 5 & 11.000 & 10.965 & 10.684 & 10.575 & 9.650 & 9.365 & 9.049 & 5.493 & 3.769 \\
 & 7 & 11.000 & 10.996 & 10.778 & 10.675 & 9.822 & 9.504 & 9.142 & 5.335 & 3.406 \\
 & 10 & 11.000 & 10.999 & 10.849 & 10.755 & 9.968 & 9.610 & 9.186 & 5.074 & 2.973 \\
 & 15 & 11.000 & 11.000 & 10.889 & 10.793 & 10.000 & 9.558 & 6.575 & 4.489 & 2.394 \\
\addlinespace[0.2em]
0.0 & 1 & 10.926 & 10.708 & 10.188 & 10.098 & 10.000 & 9.896 & 9.785 & 8.261 & 5.145 \\
 & 5 & 11.000 & 10.996 & 10.784 & 10.685 & 9.848 & 9.538 & 9.183 & 5.393 & 3.456 \\
 & 7 & 11.000 & 10.998 & 10.850 & 10.757 & 9.979 & 9.630 & 9.217 & 5.137 & 3.032 \\
 & 10 & 11.000 & 10.999 & 10.883 & 10.788 & 10.000 & 9.577 & 6.636 & 4.612 & 2.500 \\
 & 15 & 11.000 & 11.000 & 10.976 & 10.816 & 10.000 & 9.465 & 6.280 & 3.934 & 1.983 \\
\addlinespace[0.2em]
0.4 & 1 & 10.976 & 10.853 & 10.448 & 10.358 & 10.257 & 9.228 & 9.034 & 6.077 & 4.860 \\
 & 5 & 11.000 & 10.999 & 10.867 & 10.775 & 10.000 & 9.624 & 6.790 & 4.931 & 2.802 \\
 & 7 & 11.000 & 11.000 & 10.894 & 10.798 & 10.000 & 9.542 & 6.516 & 4.379 & 2.305 \\
 & 10 & 11.000 & 11.000 & 10.979 & 10.823 & 10.000 & 7.106 & 6.186 & 3.772 & 1.879 \\
 & 15 & 11.000 & 11.000 & 10.987 & 10.852 & 10.000 & 6.879 & 5.747 & 3.094 & 1.514 \\
\addlinespace[0.2em]
0.8 & 1 & 10.999 & 10.979 & 10.802 & 10.204 & 10.000 & 9.762 & 7.219 & 5.942 & 4.085 \\
 & 5 & 11.000 & 11.000 & 10.986 & 10.850 & 10.000 & 6.894 & 5.782 & 3.147 & 1.535 \\
 & 7 & 11.000 & 11.000 & 10.992 & 10.875 & 10.000 & 6.671 & 5.356 & 2.610 & 1.308 \\
 & 10 & 11.000 & 11.000 & 10.996 & 10.902 & 10.000 & 6.393 & 4.849 & 2.126 & 1.147 \\
 & 15 & 11.000 & 11.000 & 10.998 & 10.930 & 10.000 & 6.016 & 4.212 & 1.695 & 1.047 \\
\addlinespace[0.2em]
\end{longtable}
\endgroup

\begingroup
\setlength{\tabcolsep}{3pt}
\begin{longtable}{c c r r r r r r r r r}
\caption{$TARL_1$ values for the EWMA-RZ$^+$ chart with smoothing constant $\lambda = 0.2$, under equal coefficients of variation ($\gamma_X = \gamma_Y$), horizon $I = 10$, target $TARL_0 = 10$, $z_0 = 1$, across nine shift levels $\tau$.}\label{tab:tarl1_equal_T10_lam02}\\
\toprule
$\rho_0$ & $n$ & $\tau=0.9$ & $\tau=0.95$ & $\tau=0.98$ & $\tau=0.99$ & $\tau=1$ & $\tau=1.01$ & $\tau=1.02$ & $\tau=1.05$ & $\tau=1.1$ \\
\midrule
\endfirsthead
\multicolumn{11}{l}{\textit{Table \ref{tab:tarl1_equal_T10_lam02} (continued)}}\\
\toprule
$\rho_0$ & $n$ & $\tau=0.9$ & $\tau=0.95$ & $\tau=0.98$ & $\tau=0.99$ & $\tau=1$ & $\tau=1.01$ & $\tau=1.02$ & $\tau=1.05$ & $\tau=1.1$ \\
\midrule
\endhead
\bottomrule
\endlastfoot
\multicolumn{11}{l}{\textit{$\gamma_X = 0.01$, $\gamma_Y = 0.01$}} \\
-0.8 & 1 & 11.000 & 11.000 & 11.000 & 10.983 & 10.000 & 4.148 & 2.126 & 1.055 & 1.000 \\
 & 5 & 11.000 & 11.000 & 11.000 & 11.000 & 10.000 & 1.861 & 1.090 & 1.000 & 1.000 \\
 & 7 & 11.000 & 11.000 & 11.000 & 11.000 & 10.000 & 1.568 & 1.035 & 1.000 & 1.000 \\
 & 10 & 11.000 & 11.000 & 11.000 & 11.000 & 10.000 & 1.340 & 1.009 & 1.000 & 1.000 \\
 & 15 & 11.000 & 11.000 & 11.000 & 11.000 & 9.999 & 1.166 & 1.001 & 1.000 & 1.000 \\
\addlinespace[0.2em]
-0.4 & 1 & 11.000 & 11.000 & 11.000 & 10.990 & 10.000 & 3.754 & 1.852 & 1.025 & 1.000 \\
 & 5 & 11.000 & 11.000 & 11.000 & 11.000 & 10.000 & 1.635 & 1.045 & 1.000 & 1.000 \\
 & 7 & 11.000 & 11.000 & 11.000 & 11.000 & 10.000 & 1.399 & 1.014 & 1.000 & 1.000 \\
 & 10 & 11.000 & 11.000 & 11.000 & 11.000 & 10.000 & 1.222 & 1.002 & 1.000 & 1.000 \\
 & 15 & 11.000 & 11.000 & 11.000 & 11.000 & 9.999 & 1.095 & 1.000 & 1.000 & 1.000 \\
\addlinespace[0.2em]
0.0 & 1 & 11.000 & 11.000 & 11.000 & 10.998 & 10.000 & 3.220 & 1.560 & 1.006 & 1.000 \\
 & 5 & 11.000 & 11.000 & 11.000 & 11.000 & 10.000 & 1.399 & 1.014 & 1.000 & 1.000 \\
 & 7 & 11.000 & 11.000 & 11.000 & 11.000 & 10.000 & 1.231 & 1.003 & 1.000 & 1.000 \\
 & 10 & 11.000 & 11.000 & 11.000 & 11.000 & 9.999 & 1.112 & 1.000 & 1.000 & 1.000 \\
 & 15 & 11.000 & 11.000 & 11.000 & 11.000 & 9.998 & 1.037 & 1.000 & 1.000 & 1.000 \\
\addlinespace[0.2em]
0.4 & 1 & 11.000 & 11.000 & 11.000 & 11.000 & 10.000 & 2.463 & 1.258 & 1.000 & 1.000 \\
 & 5 & 11.000 & 11.000 & 11.000 & 11.000 & 9.999 & 1.166 & 1.001 & 1.000 & 1.000 \\
 & 7 & 11.000 & 11.000 & 11.000 & 11.000 & 9.999 & 1.077 & 1.000 & 1.000 & 1.000 \\
 & 10 & 11.000 & 11.000 & 11.000 & 11.000 & 9.998 & 1.026 & 1.000 & 1.000 & 1.000 \\
 & 15 & 11.000 & 11.000 & 11.000 & 11.000 & 10.003 & 1.004 & 1.000 & 1.000 & 1.000 \\
\addlinespace[0.2em]
0.8 & 1 & 11.000 & 11.000 & 11.000 & 11.000 & 10.000 & 1.400 & 1.014 & 1.000 & 1.000 \\
 & 5 & 11.000 & 11.000 & 11.000 & 11.000 & 10.003 & 1.004 & 1.000 & 1.000 & 1.000 \\
 & 7 & 11.000 & 11.000 & 11.000 & 11.000 & 10.003 & 1.001 & 1.000 & 1.000 & 1.000 \\
 & 10 & 11.000 & 11.000 & 11.000 & 11.000 & 10.003 & 1.000 & 1.000 & 1.000 & 1.000 \\
 & 15 & 11.000 & 11.000 & 11.000 & 11.000 & 10.003 & 1.000 & 1.000 & 1.000 & 1.000 \\
\addlinespace[0.2em]
\midrule
\multicolumn{11}{l}{\textit{$\gamma_X = 0.2$, $\gamma_Y = 0.2$}} \\
-0.8 & 1 & 10.749 & 10.436 & 10.286 & 10.230 & 10.000 & 9.928 & 9.853 & 9.328 & 8.375 \\
 & 5 & 10.987 & 10.871 & 10.556 & 10.187 & 10.000 & 9.787 & 8.984 & 6.985 & 3.882 \\
 & 7 & 10.997 & 10.910 & 10.407 & 10.223 & 10.000 & 9.204 & 8.843 & 6.542 & 3.272 \\
 & 10 & 10.999 & 10.955 & 10.549 & 10.361 & 10.121 & 9.300 & 8.870 & 4.760 & 2.818 \\
 & 15 & 11.000 & 10.984 & 10.817 & 10.526 & 9.875 & 9.445 & 7.951 & 4.411 & 2.360 \\
\addlinespace[0.2em]
-0.4 & 1 & 10.768 & 10.524 & 10.223 & 10.153 & 10.079 & 9.764 & 9.667 & 8.965 & 7.681 \\
 & 5 & 10.993 & 10.899 & 10.389 & 10.212 & 10.000 & 9.226 & 8.886 & 6.670 & 3.429 \\
 & 7 & 10.999 & 10.943 & 10.508 & 10.320 & 10.086 & 9.275 & 8.868 & 4.860 & 2.955 \\
 & 10 & 11.000 & 10.976 & 10.654 & 10.475 & 9.818 & 9.413 & 8.937 & 4.581 & 2.546 \\
 & 15 & 11.000 & 10.993 & 10.867 & 10.610 & 9.983 & 9.516 & 7.949 & 4.148 & 2.096 \\
\addlinespace[0.2em]
0.0 & 1 & 10.813 & 10.549 & 10.183 & 10.095 & 10.000 & 9.898 & 9.449 & 9.022 & 7.547 \\
 & 5 & 10.999 & 10.944 & 10.520 & 10.339 & 10.111 & 9.316 & 8.918 & 6.431 & 2.994 \\
 & 7 & 11.000 & 10.975 & 10.654 & 10.479 & 9.831 & 9.437 & 8.972 & 4.642 & 2.594 \\
 & 10 & 11.000 & 10.991 & 10.858 & 10.597 & 9.974 & 9.525 & 7.993 & 4.260 & 2.183 \\
 & 15 & 11.000 & 10.998 & 10.900 & 10.656 & 10.000 & 9.435 & 7.681 & 3.572 & 1.734 \\
\addlinespace[0.2em]
0.4 & 1 & 10.923 & 10.654 & 10.450 & 10.124 & 10.000 & 9.864 & 9.715 & 8.615 & 6.423 \\
 & 5 & 11.000 & 10.985 & 10.723 & 10.561 & 9.951 & 9.547 & 9.058 & 4.555 & 2.436 \\
 & 7 & 11.000 & 10.994 & 10.876 & 10.625 & 10.000 & 9.519 & 7.923 & 4.034 & 2.011 \\
 & 10 & 11.000 & 10.999 & 10.908 & 10.668 & 10.000 & 9.401 & 7.583 & 3.392 & 1.642 \\
 & 15 & 11.000 & 11.000 & 10.940 & 10.719 & 10.000 & 8.277 & 7.111 & 2.713 & 1.354 \\
\addlinespace[0.2em]
0.8 & 1 & 10.992 & 10.907 & 10.627 & 10.213 & 10.000 & 9.746 & 9.449 & 7.287 & 3.535 \\
 & 5 & 11.000 & 11.000 & 10.936 & 10.713 & 10.000 & 8.294 & 7.144 & 2.733 & 1.362 \\
 & 7 & 11.000 & 11.000 & 10.959 & 10.759 & 10.000 & 8.072 & 4.948 & 2.243 & 1.199 \\
 & 10 & 11.000 & 11.000 & 10.977 & 10.807 & 10.000 & 7.788 & 4.408 & 1.829 & 1.090 \\
 & 15 & 11.000 & 11.000 & 10.990 & 10.859 & 10.000 & 7.390 & 3.751 & 1.486 & 1.027 \\
\addlinespace[0.2em]
\end{longtable}
\endgroup

\begingroup
\setlength{\tabcolsep}{3pt}
\begin{longtable}{c c r r r r r r r r r}
\caption{$TARL_1$ values for the EWMA-RZ$^+$ chart with smoothing constant $\lambda = 0.1$, under unequal coefficients of variation ($\gamma_X \neq \gamma_Y$), horizon $I = 10$, target $TARL_0 = 10$, $z_0 = 1$, across nine shift levels $\tau$.}\label{tab:tarl1_unequal_T10_lam01}\\
\toprule
$\rho_0$ & $n$ & $\tau=0.9$ & $\tau=0.95$ & $\tau=0.98$ & $\tau=0.99$ & $\tau=1$ & $\tau=1.01$ & $\tau=1.02$ & $\tau=1.05$ & $\tau=1.1$ \\
\midrule
\endfirsthead
\multicolumn{11}{l}{\textit{Table \ref{tab:tarl1_unequal_T10_lam01} (continued)}}\\
\toprule
$\rho_0$ & $n$ & $\tau=0.9$ & $\tau=0.95$ & $\tau=0.98$ & $\tau=0.99$ & $\tau=1$ & $\tau=1.01$ & $\tau=1.02$ & $\tau=1.05$ & $\tau=1.1$ \\
\midrule
\endhead
\bottomrule
\endlastfoot
\multicolumn{11}{l}{\textit{$\gamma_X = 0.01$, $\gamma_Y = 0.2$}} \\
-0.8 & 1 & 10.959 & 10.795 & 10.504 & 10.116 & 10.000 & 9.872 & 9.078 & 7.269 & 4.179 \\
 & 5 & 11.000 & 10.998 & 10.881 & 10.517 & 10.273 & 8.878 & 5.922 & 4.081 & 2.275 \\
 & 7 & 11.000 & 11.000 & 10.934 & 10.655 & 9.611 & 9.039 & 5.928 & 3.798 & 1.984 \\
 & 10 & 11.000 & 11.000 & 10.967 & 10.761 & 9.803 & 9.152 & 5.857 & 3.421 & 1.698 \\
 & 15 & 11.000 & 11.000 & 10.986 & 10.845 & 9.973 & 9.212 & 5.670 & 2.929 & 1.424 \\
\addlinespace[0.2em]
-0.4 & 1 & 10.960 & 10.798 & 10.506 & 10.118 & 10.000 & 9.870 & 9.071 & 7.244 & 4.136 \\
 & 5 & 11.000 & 10.998 & 10.883 & 10.521 & 10.273 & 8.869 & 5.898 & 4.031 & 2.231 \\
 & 7 & 11.000 & 11.000 & 10.936 & 10.658 & 9.612 & 9.028 & 5.898 & 3.742 & 1.945 \\
 & 10 & 11.000 & 11.000 & 10.968 & 10.764 & 9.803 & 9.139 & 5.822 & 3.363 & 1.665 \\
 & 15 & 11.000 & 11.000 & 10.987 & 10.848 & 9.973 & 9.196 & 5.626 & 2.871 & 1.399 \\
\addlinespace[0.2em]
0.0 & 1 & 10.961 & 10.800 & 10.509 & 10.120 & 10.000 & 9.867 & 9.064 & 7.219 & 4.091 \\
 & 5 & 11.000 & 10.999 & 10.886 & 10.525 & 10.274 & 8.859 & 5.872 & 3.979 & 2.187 \\
 & 7 & 11.000 & 11.000 & 10.938 & 10.662 & 9.612 & 9.016 & 5.867 & 3.685 & 1.905 \\
 & 10 & 11.000 & 11.000 & 10.969 & 10.767 & 9.803 & 9.125 & 5.784 & 3.301 & 1.632 \\
 & 15 & 11.000 & 11.000 & 10.987 & 10.851 & 9.973 & 9.178 & 5.580 & 2.811 & 1.374 \\
\addlinespace[0.2em]
0.4 & 1 & 10.962 & 10.803 & 10.511 & 10.122 & 10.000 & 9.865 & 9.056 & 7.193 & 4.044 \\
 & 5 & 11.000 & 10.999 & 10.888 & 10.529 & 10.274 & 8.847 & 5.844 & 3.923 & 2.141 \\
 & 7 & 11.000 & 11.000 & 10.940 & 10.666 & 9.612 & 9.003 & 5.833 & 3.624 & 1.864 \\
 & 10 & 11.000 & 11.000 & 10.970 & 10.771 & 9.803 & 9.109 & 5.743 & 3.238 & 1.598 \\
 & 15 & 11.000 & 11.000 & 10.988 & 10.854 & 9.973 & 9.159 & 5.530 & 2.749 & 1.348 \\
\addlinespace[0.2em]
0.8 & 1 & 10.963 & 10.806 & 10.514 & 10.125 & 10.000 & 9.862 & 9.048 & 7.166 & 3.996 \\
 & 5 & 11.000 & 10.999 & 10.891 & 10.533 & 10.273 & 8.833 & 5.812 & 3.863 & 2.094 \\
 & 7 & 11.000 & 11.000 & 10.941 & 10.670 & 9.611 & 8.988 & 5.796 & 3.559 & 1.823 \\
 & 10 & 11.000 & 11.000 & 10.972 & 10.774 & 9.803 & 9.092 & 5.699 & 3.171 & 1.564 \\
 & 15 & 11.000 & 11.000 & 10.989 & 10.857 & 9.973 & 6.760 & 5.477 & 2.685 & 1.323 \\
\addlinespace[0.2em]
\midrule
\multicolumn{11}{l}{\textit{$\gamma_X = 0.2$, $\gamma_Y = 0.01$}} \\
-0.8 & 1 & 10.998 & 10.912 & 10.311 & 10.164 & 10.000 & 7.642 & 7.333 & 6.407 & 4.997 \\
 & 5 & 11.000 & 10.998 & 10.889 & 10.331 & 10.000 & 7.268 & 6.571 & 4.620 & 2.584 \\
 & 7 & 11.000 & 10.999 & 10.912 & 10.381 & 10.000 & 7.142 & 6.315 & 4.106 & 2.151 \\
 & 10 & 11.000 & 11.000 & 10.935 & 10.439 & 10.000 & 6.982 & 5.994 & 3.537 & 1.772 \\
 & 15 & 11.000 & 11.000 & 10.958 & 10.513 & 10.000 & 6.760 & 5.559 & 2.901 & 1.447 \\
\addlinespace[0.2em]
-0.4 & 1 & 10.998 & 10.915 & 10.316 & 10.167 & 10.000 & 7.634 & 7.319 & 6.376 & 4.947 \\
 & 5 & 11.000 & 10.998 & 10.891 & 10.337 & 10.000 & 7.254 & 6.543 & 4.563 & 2.533 \\
 & 7 & 11.000 & 10.999 & 10.914 & 10.387 & 10.000 & 7.125 & 6.283 & 4.046 & 2.107 \\
 & 10 & 11.000 & 11.000 & 10.937 & 10.446 & 10.000 & 6.962 & 5.956 & 3.476 & 1.738 \\
 & 15 & 11.000 & 11.000 & 10.960 & 10.520 & 10.000 & 6.736 & 5.513 & 2.845 & 1.422 \\
\addlinespace[0.2em]
0.0 & 1 & 10.998 & 10.918 & 10.321 & 10.170 & 10.000 & 7.626 & 7.305 & 6.344 & 4.895 \\
 & 5 & 11.000 & 10.998 & 10.894 & 10.343 & 10.000 & 7.239 & 6.513 & 4.504 & 2.481 \\
 & 7 & 11.000 & 10.999 & 10.917 & 10.394 & 10.000 & 7.108 & 6.248 & 3.983 & 2.063 \\
 & 10 & 11.000 & 11.000 & 10.940 & 10.453 & 10.000 & 6.941 & 5.915 & 3.413 & 1.704 \\
 & 15 & 11.000 & 11.000 & 10.962 & 10.528 & 10.000 & 6.710 & 5.465 & 2.787 & 1.398 \\
\addlinespace[0.2em]
0.4 & 1 & 10.998 & 10.924 & 10.348 & 10.197 & 10.027 & 7.667 & 7.340 & 6.360 & 4.883 \\
 & 5 & 11.000 & 10.999 & 10.601 & 10.349 & 10.000 & 7.223 & 6.482 & 4.442 & 2.428 \\
 & 7 & 11.000 & 10.999 & 10.920 & 10.401 & 10.000 & 7.089 & 6.211 & 3.918 & 2.019 \\
 & 10 & 11.000 & 11.000 & 10.942 & 10.461 & 10.000 & 6.919 & 5.872 & 3.348 & 1.669 \\
 & 15 & 11.000 & 11.000 & 10.964 & 10.536 & 10.000 & 6.683 & 5.414 & 2.727 & 1.373 \\
\addlinespace[0.2em]
0.8 & 1 & 10.999 & 10.935 & 10.402 & 10.256 & 10.090 & 7.778 & 7.446 & 6.446 & 4.930 \\
 & 5 & 11.000 & 10.999 & 10.623 & 10.374 & 10.024 & 7.252 & 6.495 & 4.414 & 2.389 \\
 & 7 & 11.000 & 11.000 & 10.925 & 10.419 & 10.016 & 7.099 & 6.202 & 3.871 & 1.981 \\
 & 10 & 11.000 & 11.000 & 10.946 & 10.474 & 10.009 & 6.912 & 5.842 & 3.291 & 1.636 \\
 & 15 & 11.000 & 11.000 & 10.966 & 10.546 & 10.002 & 6.657 & 5.363 & 2.667 & 1.349 \\
\addlinespace[0.2em]
\end{longtable}
\endgroup

\begingroup
\setlength{\tabcolsep}{3pt}
\begin{longtable}{c c r r r r r r r r r}
\caption{$TARL_1$ values for the EWMA-RZ$^+$ chart with smoothing constant $\lambda = 0.2$, under unequal coefficients of variation ($\gamma_X \neq \gamma_Y$), horizon $I = 10$, target $TARL_0 = 10$, $z_0 = 1$, across nine shift levels $\tau$.}\label{tab:tarl1_unequal_T10_lam02}\\
\toprule
$\rho_0$ & $n$ & $\tau=0.9$ & $\tau=0.95$ & $\tau=0.98$ & $\tau=0.99$ & $\tau=1$ & $\tau=1.01$ & $\tau=1.02$ & $\tau=1.05$ & $\tau=1.1$ \\
\midrule
\endfirsthead
\multicolumn{11}{l}{\textit{Table \ref{tab:tarl1_unequal_T10_lam02} (continued)}}\\
\toprule
$\rho_0$ & $n$ & $\tau=0.9$ & $\tau=0.95$ & $\tau=0.98$ & $\tau=0.99$ & $\tau=1$ & $\tau=1.01$ & $\tau=1.02$ & $\tau=1.05$ & $\tau=1.1$ \\
\midrule
\endhead
\bottomrule
\endlastfoot
\multicolumn{11}{l}{\textit{$\gamma_X = 0.01$, $\gamma_Y = 0.2$}} \\
-0.8 & 1 & 10.847 & 10.577 & 10.266 & 10.174 & 10.073 & 9.766 & 9.630 & 8.833 & 6.674 \\
 & 5 & 11.000 & 10.974 & 10.691 & 10.522 & 10.000 & 9.067 & 8.510 & 4.214 & 2.211 \\
 & 7 & 11.000 & 10.988 & 10.752 & 10.361 & 10.000 & 8.937 & 7.221 & 3.561 & 1.786 \\
 & 10 & 11.000 & 10.996 & 10.830 & 10.471 & 10.068 & 8.886 & 6.949 & 3.009 & 1.493 \\
 & 15 & 11.000 & 11.000 & 10.916 & 10.643 & 9.824 & 8.976 & 6.790 & 2.537 & 1.282 \\
\addlinespace[0.2em]
-0.4 & 1 & 10.850 & 10.582 & 10.270 & 10.177 & 10.074 & 9.766 & 9.627 & 8.820 & 6.633 \\
 & 5 & 11.000 & 10.975 & 10.696 & 10.525 & 10.000 & 9.058 & 8.488 & 4.160 & 2.164 \\
 & 7 & 11.000 & 10.989 & 10.757 & 10.366 & 10.000 & 8.924 & 7.191 & 3.505 & 1.751 \\
 & 10 & 11.000 & 10.997 & 10.835 & 10.477 & 10.068 & 8.872 & 6.912 & 2.954 & 1.467 \\
 & 15 & 11.000 & 11.000 & 10.919 & 10.649 & 9.824 & 8.958 & 6.743 & 2.484 & 1.263 \\
\addlinespace[0.2em]
0.0 & 1 & 10.854 & 10.587 & 10.273 & 10.179 & 10.075 & 9.765 & 9.623 & 8.804 & 6.587 \\
 & 5 & 11.000 & 10.976 & 10.701 & 10.529 & 10.000 & 9.048 & 8.465 & 4.103 & 2.117 \\
 & 7 & 11.000 & 10.990 & 10.762 & 10.372 & 10.000 & 8.912 & 7.158 & 3.446 & 1.716 \\
 & 10 & 11.000 & 10.998 & 10.840 & 10.484 & 10.069 & 8.857 & 6.873 & 2.897 & 1.442 \\
 & 15 & 11.000 & 11.000 & 10.923 & 10.655 & 9.825 & 8.939 & 6.694 & 2.431 & 1.245 \\
\addlinespace[0.2em]
0.4 & 1 & 10.879 & 10.591 & 10.276 & 10.181 & 10.075 & 9.762 & 9.618 & 8.785 & 6.536 \\
 & 5 & 11.000 & 10.978 & 10.707 & 10.533 & 10.000 & 9.037 & 8.441 & 4.043 & 2.069 \\
 & 7 & 11.000 & 10.991 & 10.768 & 10.378 & 10.000 & 8.898 & 7.124 & 3.386 & 1.680 \\
 & 10 & 11.000 & 10.999 & 10.845 & 10.490 & 10.069 & 8.840 & 6.830 & 2.839 & 1.416 \\
 & 15 & 11.000 & 11.000 & 10.926 & 10.660 & 9.825 & 8.918 & 6.641 & 2.376 & 1.227 \\
\addlinespace[0.2em]
0.8 & 1 & 10.881 & 10.595 & 10.278 & 10.181 & 10.073 & 9.757 & 9.610 & 8.762 & 6.480 \\
 & 5 & 11.000 & 10.979 & 10.712 & 10.537 & 10.000 & 9.026 & 8.414 & 3.981 & 2.021 \\
 & 7 & 11.000 & 10.991 & 10.774 & 10.384 & 10.000 & 8.884 & 7.087 & 3.323 & 1.644 \\
 & 10 & 11.000 & 10.999 & 10.850 & 10.496 & 10.068 & 8.820 & 6.783 & 2.777 & 1.390 \\
 & 15 & 11.000 & 11.000 & 10.929 & 10.666 & 9.824 & 8.895 & 6.583 & 2.319 & 1.209 \\
\addlinespace[0.2em]
\midrule
\multicolumn{11}{l}{\textit{$\gamma_X = 0.2$, $\gamma_Y = 0.01$}} \\
-0.8 & 1 & 10.989 & 10.839 & 10.341 & 10.182 & 10.000 & 9.796 & 8.671 & 7.716 & 4.385 \\
 & 5 & 11.000 & 10.989 & 10.791 & 10.356 & 10.000 & 8.621 & 7.925 & 4.047 & 2.154 \\
 & 7 & 11.000 & 10.995 & 10.832 & 10.408 & 10.000 & 8.502 & 7.661 & 3.545 & 1.808 \\
 & 10 & 11.000 & 10.999 & 10.874 & 10.468 & 10.000 & 8.347 & 5.481 & 3.011 & 1.523 \\
 & 15 & 11.000 & 11.000 & 10.916 & 10.542 & 10.000 & 8.129 & 5.027 & 2.444 & 1.290 \\
\addlinespace[0.2em]
-0.4 & 1 & 10.990 & 10.845 & 10.347 & 10.186 & 10.000 & 8.955 & 8.654 & 7.676 & 4.308 \\
 & 5 & 11.000 & 10.990 & 10.797 & 10.363 & 10.000 & 8.606 & 7.893 & 3.980 & 2.108 \\
 & 7 & 11.000 & 10.996 & 10.837 & 10.415 & 10.000 & 8.484 & 7.624 & 3.480 & 1.772 \\
 & 10 & 11.000 & 10.999 & 10.878 & 10.475 & 10.000 & 8.327 & 5.432 & 2.950 & 1.497 \\
 & 15 & 11.000 & 11.000 & 10.920 & 10.550 & 10.000 & 8.103 & 4.973 & 2.393 & 1.272 \\
\addlinespace[0.2em]
0.0 & 1 & 10.991 & 10.851 & 10.354 & 10.190 & 10.000 & 8.946 & 8.637 & 7.634 & 4.230 \\
 & 5 & 11.000 & 10.991 & 10.802 & 10.369 & 10.000 & 8.590 & 7.860 & 3.911 & 2.062 \\
 & 7 & 11.000 & 10.996 & 10.843 & 10.422 & 10.000 & 8.466 & 7.584 & 3.413 & 1.737 \\
 & 10 & 11.000 & 11.000 & 10.883 & 10.483 & 10.000 & 8.305 & 5.381 & 2.888 & 1.471 \\
 & 15 & 11.000 & 11.000 & 10.924 & 10.558 & 10.000 & 8.076 & 4.917 & 2.340 & 1.255 \\
\addlinespace[0.2em]
0.4 & 1 & 10.992 & 10.857 & 10.362 & 10.194 & 10.000 & 8.936 & 8.619 & 5.715 & 4.150 \\
 & 5 & 11.000 & 10.992 & 10.808 & 10.376 & 10.000 & 8.573 & 7.824 & 3.840 & 2.015 \\
 & 7 & 11.000 & 10.997 & 10.848 & 10.429 & 10.000 & 8.446 & 5.679 & 3.343 & 1.701 \\
 & 10 & 11.000 & 11.000 & 10.888 & 10.491 & 10.000 & 8.281 & 5.327 & 2.824 & 1.445 \\
 & 15 & 11.000 & 11.000 & 10.928 & 10.567 & 10.000 & 8.047 & 4.858 & 2.287 & 1.238 \\
\addlinespace[0.2em]
0.8 & 1 & 10.994 & 10.867 & 10.382 & 10.213 & 10.017 & 8.953 & 8.628 & 5.679 & 4.092 \\
 & 5 & 11.000 & 10.993 & 10.814 & 10.384 & 10.000 & 8.555 & 7.786 & 3.765 & 1.968 \\
 & 7 & 11.000 & 10.997 & 10.854 & 10.437 & 10.000 & 8.425 & 5.627 & 3.272 & 1.665 \\
 & 10 & 11.000 & 11.000 & 10.894 & 10.500 & 10.000 & 8.256 & 5.270 & 2.758 & 1.419 \\
 & 15 & 11.000 & 11.000 & 10.933 & 10.576 & 10.000 & 8.017 & 4.795 & 2.233 & 1.221 \\
\addlinespace[0.2em]
\end{longtable}
\endgroup

\begingroup
\setlength{\tabcolsep}{3pt}
\begin{longtable}{c c r r r r r r r r r}
\caption{$TARL_1$ values for the EWMA-RZ$^+$ chart with smoothing constant $\lambda = 0.1$, under equal coefficients of variation ($\gamma_X = \gamma_Y$), horizon $I = 30$, target $TARL_0 = 30$, $z_0 = 1$, across nine shift levels $\tau$.}\label{tab:tarl1_equal_T30_lam01}\\
\toprule
$\rho_0$ & $n$ & $\tau=0.9$ & $\tau=0.95$ & $\tau=0.98$ & $\tau=0.99$ & $\tau=1$ & $\tau=1.01$ & $\tau=1.02$ & $\tau=1.05$ & $\tau=1.1$ \\
\midrule
\endfirsthead
\multicolumn{11}{l}{\textit{Table \ref{tab:tarl1_equal_T30_lam01} (continued)}}\\
\toprule
$\rho_0$ & $n$ & $\tau=0.9$ & $\tau=0.95$ & $\tau=0.98$ & $\tau=0.99$ & $\tau=1$ & $\tau=1.01$ & $\tau=1.02$ & $\tau=1.05$ & $\tau=1.1$ \\
\midrule
\endhead
\bottomrule
\endlastfoot
\multicolumn{11}{l}{\textit{$\gamma_X = 0.01$, $\gamma_Y = 0.01$}} \\
-0.8 & 1 & 31.000 & 31.000 & 31.000 & 30.999 & 30.000 & 11.197 & 3.903 & 1.178 & 1.000 \\
 & 5 & 31.000 & 31.000 & 31.000 & 31.000 & 29.999 & 3.235 & 1.275 & 1.000 & 1.000 \\
 & 7 & 31.000 & 31.000 & 31.000 & 31.000 & 29.999 & 2.479 & 1.121 & 1.000 & 1.000 \\
 & 10 & 31.000 & 31.000 & 31.000 & 31.000 & 29.999 & 1.908 & 1.037 & 1.000 & 1.000 \\
 & 15 & 31.000 & 31.000 & 31.000 & 31.000 & 29.999 & 1.472 & 1.005 & 1.000 & 1.000 \\
\addlinespace[0.2em]
-0.4 & 1 & 31.000 & 31.000 & 31.000 & 31.000 & 30.000 & 9.591 & 3.177 & 1.089 & 1.000 \\
 & 5 & 31.000 & 31.000 & 31.000 & 31.000 & 29.999 & 2.648 & 1.152 & 1.000 & 1.000 \\
 & 7 & 31.000 & 31.000 & 31.000 & 31.000 & 29.999 & 2.056 & 1.054 & 1.000 & 1.000 \\
 & 10 & 31.000 & 31.000 & 31.000 & 31.000 & 29.999 & 1.615 & 1.012 & 1.000 & 1.000 \\
 & 15 & 31.000 & 31.000 & 31.000 & 31.000 & 29.999 & 1.290 & 1.001 & 1.000 & 1.000 \\
\addlinespace[0.2em]
0.0 & 1 & 31.000 & 31.000 & 31.000 & 31.000 & 29.999 & 7.531 & 2.440 & 1.026 & 1.000 \\
 & 5 & 31.000 & 31.000 & 31.000 & 31.000 & 29.999 & 2.056 & 1.054 & 1.000 & 1.000 \\
 & 7 & 31.000 & 31.000 & 31.000 & 31.000 & 29.999 & 1.636 & 1.013 & 1.000 & 1.000 \\
 & 10 & 31.000 & 31.000 & 31.000 & 31.000 & 29.999 & 1.334 & 1.001 & 1.000 & 1.000 \\
 & 15 & 31.000 & 31.000 & 31.000 & 31.000 & 29.999 & 1.129 & 1.000 & 1.000 & 1.000 \\
\addlinespace[0.2em]
0.4 & 1 & 31.000 & 31.000 & 31.000 & 31.000 & 29.999 & 4.958 & 1.699 & 1.001 & 1.000 \\
 & 5 & 31.000 & 31.000 & 31.000 & 31.000 & 29.999 & 1.472 & 1.005 & 1.000 & 1.000 \\
 & 7 & 31.000 & 31.000 & 31.000 & 31.000 & 29.999 & 1.241 & 1.000 & 1.000 & 1.000 \\
 & 10 & 31.000 & 31.000 & 31.000 & 31.000 & 29.999 & 1.095 & 1.000 & 1.000 & 1.000 \\
 & 15 & 31.000 & 31.000 & 31.000 & 31.000 & 29.999 & 1.021 & 1.000 & 1.000 & 1.000 \\
\addlinespace[0.2em]
0.8 & 1 & 31.000 & 31.000 & 31.000 & 31.000 & 29.999 & 2.057 & 1.054 & 1.000 & 1.000 \\
 & 5 & 31.000 & 31.000 & 31.000 & 31.000 & 29.999 & 1.021 & 1.000 & 1.000 & 1.000 \\
 & 7 & 31.000 & 31.000 & 31.000 & 31.000 & 29.999 & 1.003 & 1.000 & 1.000 & 1.000 \\
 & 10 & 31.000 & 31.000 & 31.000 & 31.000 & 29.999 & 1.000 & 1.000 & 1.000 & 1.000 \\
 & 15 & 31.000 & 31.000 & 31.000 & 31.000 & 29.999 & 1.000 & 1.000 & 1.000 & 1.000 \\
\addlinespace[0.2em]
\midrule
\multicolumn{11}{l}{\textit{$\gamma_X = 0.2$, $\gamma_Y = 0.2$}} \\
-0.8 & 1 & 30.842 & 30.600 & 30.241 & 30.128 & 30.000 & 29.856 & 29.390 & 28.685 & 26.240 \\
 & 5 & 31.000 & 30.974 & 30.747 & 30.581 & 30.000 & 29.521 & 28.883 & 24.207 & 8.730 \\
 & 7 & 31.000 & 30.987 & 30.794 & 30.401 & 30.000 & 29.409 & 27.457 & 19.683 & 6.804 \\
 & 10 & 31.000 & 30.994 & 30.841 & 30.464 & 30.000 & 28.342 & 26.947 & 17.633 & 5.151 \\
 & 15 & 31.000 & 30.998 & 30.890 & 30.539 & 30.000 & 28.063 & 26.187 & 9.895 & 3.763 \\
\addlinespace[0.2em]
-0.4 & 1 & 30.905 & 30.658 & 30.432 & 30.146 & 30.000 & 29.832 & 29.639 & 28.368 & 25.116 \\
 & 5 & 31.000 & 30.983 & 30.780 & 30.385 & 30.000 & 29.443 & 27.577 & 20.191 & 7.324 \\
 & 7 & 31.000 & 30.992 & 30.826 & 30.443 & 30.000 & 28.409 & 27.125 & 18.323 & 5.635 \\
 & 10 & 31.000 & 30.997 & 30.871 & 30.508 & 30.000 & 28.182 & 26.515 & 10.804 & 4.256 \\
 & 15 & 31.000 & 31.000 & 30.926 & 30.630 & 30.092 & 28.069 & 23.290 & 8.786 & 3.190 \\
\addlinespace[0.2em]
0.0 & 1 & 30.943 & 30.799 & 30.502 & 30.388 & 30.000 & 29.794 & 29.552 & 27.818 & 23.120 \\
 & 5 & 31.000 & 30.992 & 30.823 & 30.440 & 30.000 & 28.420 & 27.150 & 18.419 & 5.707 \\
 & 7 & 31.000 & 30.997 & 30.867 & 30.502 & 30.000 & 28.205 & 26.576 & 10.983 & 4.365 \\
 & 10 & 31.000 & 30.999 & 30.917 & 30.608 & 30.075 & 28.096 & 23.451 & 9.171 & 3.370 \\
 & 15 & 31.000 & 31.000 & 30.965 & 30.758 & 29.754 & 28.199 & 22.969 & 7.460 & 2.591 \\
\addlinespace[0.2em]
0.4 & 1 & 30.973 & 30.860 & 30.525 & 30.373 & 30.180 & 29.493 & 29.101 & 26.092 & 18.141 \\
 & 5 & 31.000 & 30.998 & 30.808 & 30.547 & 30.029 & 28.159 & 26.361 & 10.208 & 3.913 \\
 & 7 & 31.000 & 31.000 & 30.943 & 30.689 & 30.198 & 28.275 & 23.556 & 8.775 & 3.123 \\
 & 10 & 31.000 & 31.000 & 30.974 & 30.796 & 29.854 & 28.291 & 22.961 & 7.149 & 2.455 \\
 & 15 & 31.000 & 31.000 & 30.990 & 30.874 & 30.000 & 28.109 & 21.732 & 5.356 & 1.877 \\
\addlinespace[0.2em]
0.8 & 1 & 30.997 & 30.970 & 30.741 & 30.569 & 29.802 & 29.209 & 28.406 & 21.731 & 8.218 \\
 & 5 & 31.000 & 31.000 & 30.989 & 30.869 & 30.000 & 28.164 & 21.934 & 5.560 & 1.930 \\
 & 7 & 31.000 & 31.000 & 30.995 & 30.904 & 30.000 & 25.472 & 20.035 & 4.142 & 1.541 \\
 & 10 & 31.000 & 31.000 & 30.998 & 30.936 & 30.000 & 24.509 & 12.031 & 3.069 & 1.270 \\
 & 15 & 31.000 & 31.000 & 31.000 & 30.982 & 30.000 & 23.062 & 9.400 & 2.236 & 1.096 \\
\addlinespace[0.2em]
\end{longtable}
\endgroup

\begingroup
\setlength{\tabcolsep}{3pt}
\begin{longtable}{c c r r r r r r r r r}
\caption{$TARL_1$ values for the EWMA-RZ$^+$ chart with smoothing constant $\lambda = 0.2$, under equal coefficients of variation ($\gamma_X = \gamma_Y$), horizon $I = 30$, target $TARL_0 = 30$, $z_0 = 1$, across nine shift levels $\tau$.}\label{tab:tarl1_equal_T30_lam02}\\
\toprule
$\rho_0$ & $n$ & $\tau=0.9$ & $\tau=0.95$ & $\tau=0.98$ & $\tau=0.99$ & $\tau=1$ & $\tau=1.01$ & $\tau=1.02$ & $\tau=1.05$ & $\tau=1.1$ \\
\midrule
\endfirsthead
\multicolumn{11}{l}{\textit{Table \ref{tab:tarl1_equal_T30_lam02} (continued)}}\\
\toprule
$\rho_0$ & $n$ & $\tau=0.9$ & $\tau=0.95$ & $\tau=0.98$ & $\tau=0.99$ & $\tau=1$ & $\tau=1.01$ & $\tau=1.02$ & $\tau=1.05$ & $\tau=1.1$ \\
\midrule
\endhead
\bottomrule
\endlastfoot
\multicolumn{11}{l}{\textit{$\gamma_X = 0.01$, $\gamma_Y = 0.01$}} \\
-0.8 & 1 & 31.000 & 31.000 & 31.000 & 30.996 & 30.000 & 15.939 & 3.707 & 1.124 & 1.000 \\
 & 5 & 31.000 & 31.000 & 31.000 & 31.000 & 29.999 & 2.967 & 1.197 & 1.000 & 1.000 \\
 & 7 & 31.000 & 31.000 & 31.000 & 31.000 & 29.998 & 2.224 & 1.083 & 1.000 & 1.000 \\
 & 10 & 31.000 & 31.000 & 31.000 & 31.000 & 29.997 & 1.710 & 1.024 & 1.000 & 1.000 \\
 & 15 & 31.000 & 31.000 & 31.000 & 31.000 & 30.003 & 1.350 & 1.003 & 1.000 & 1.000 \\
\addlinespace[0.2em]
-0.4 & 1 & 31.000 & 31.000 & 31.000 & 30.998 & 30.000 & 13.819 & 2.927 & 1.060 & 1.000 \\
 & 5 & 31.000 & 31.000 & 31.000 & 31.000 & 29.998 & 2.384 & 1.105 & 1.000 & 1.000 \\
 & 7 & 31.000 & 31.000 & 31.000 & 31.000 & 29.998 & 1.839 & 1.036 & 1.000 & 1.000 \\
 & 10 & 31.000 & 31.000 & 31.000 & 31.000 & 30.003 & 1.465 & 1.007 & 1.000 & 1.000 \\
 & 15 & 31.000 & 31.000 & 31.000 & 31.000 & 30.003 & 1.208 & 1.000 & 1.000 & 1.000 \\
\addlinespace[0.2em]
0.0 & 1 & 31.000 & 31.000 & 31.000 & 31.000 & 30.000 & 7.895 & 2.196 & 1.017 & 1.000 \\
 & 5 & 31.000 & 31.000 & 31.000 & 31.000 & 29.998 & 1.839 & 1.036 & 1.000 & 1.000 \\
 & 7 & 31.000 & 31.000 & 31.000 & 31.000 & 30.003 & 1.482 & 1.008 & 1.000 & 1.000 \\
 & 10 & 31.000 & 31.000 & 31.000 & 31.000 & 30.003 & 1.242 & 1.001 & 1.000 & 1.000 \\
 & 15 & 31.000 & 31.000 & 31.000 & 31.000 & 30.003 & 1.089 & 1.000 & 1.000 & 1.000 \\
\addlinespace[0.2em]
0.4 & 1 & 31.000 & 31.000 & 31.000 & 31.000 & 30.000 & 4.856 & 1.536 & 1.001 & 1.000 \\
 & 5 & 31.000 & 31.000 & 31.000 & 31.000 & 30.003 & 1.350 & 1.003 & 1.000 & 1.000 \\
 & 7 & 31.000 & 31.000 & 31.000 & 31.000 & 30.003 & 1.171 & 1.000 & 1.000 & 1.000 \\
 & 10 & 31.000 & 31.000 & 31.000 & 31.000 & 30.003 & 1.065 & 1.000 & 1.000 & 1.000 \\
 & 15 & 31.000 & 31.000 & 31.000 & 31.000 & 30.003 & 1.014 & 1.000 & 1.000 & 1.000 \\
\addlinespace[0.2em]
0.8 & 1 & 31.000 & 31.000 & 31.000 & 31.000 & 29.998 & 1.840 & 1.036 & 1.000 & 1.000 \\
 & 5 & 31.000 & 31.000 & 31.000 & 31.000 & 30.003 & 1.014 & 1.000 & 1.000 & 1.000 \\
 & 7 & 31.000 & 31.000 & 31.000 & 31.000 & 30.006 & 1.002 & 1.000 & 1.000 & 1.000 \\
 & 10 & 31.000 & 31.000 & 31.000 & 31.000 & 30.006 & 1.000 & 1.000 & 1.000 & 1.000 \\
 & 15 & 31.000 & 31.000 & 31.000 & 31.000 & 30.006 & 1.000 & 1.000 & 1.000 & 1.000 \\
\addlinespace[0.2em]
\midrule
\multicolumn{11}{l}{\textit{$\gamma_X = 0.2$, $\gamma_Y = 0.2$}} \\
-0.8 & 1 & 30.703 & 30.454 & 30.191 & 30.100 & 30.000 & 29.808 & 29.679 & 29.082 & 27.723 \\
 & 5 & 30.997 & 30.919 & 30.614 & 30.415 & 30.000 & 29.576 & 28.697 & 25.856 & 15.989 \\
 & 7 & 30.999 & 30.954 & 30.678 & 30.370 & 30.000 & 29.203 & 28.391 & 23.642 & 10.935 \\
 & 10 & 31.000 & 30.981 & 30.733 & 30.500 & 29.934 & 29.239 & 28.235 & 20.821 & 6.188 \\
 & 15 & 31.000 & 30.994 & 30.821 & 30.598 & 30.000 & 29.133 & 27.202 & 18.192 & 4.169 \\
\addlinespace[0.2em]
-0.4 & 1 & 30.772 & 30.542 & 30.280 & 30.186 & 30.000 & 29.872 & 29.730 & 29.044 & 27.388 \\
 & 5 & 30.999 & 30.944 & 30.658 & 30.352 & 30.000 & 29.241 & 28.493 & 24.111 & 11.802 \\
 & 7 & 31.000 & 30.972 & 30.678 & 30.427 & 30.023 & 29.146 & 28.176 & 21.123 & 8.995 \\
 & 10 & 31.000 & 30.991 & 30.795 & 30.572 & 30.000 & 29.223 & 27.491 & 19.432 & 4.911 \\
 & 15 & 31.000 & 30.997 & 30.859 & 30.638 & 30.000 & 28.976 & 26.682 & 13.959 & 3.264 \\
\addlinespace[0.2em]
0.0 & 1 & 30.832 & 30.575 & 30.328 & 30.138 & 30.000 & 29.841 & 29.659 & 28.749 & 26.388 \\
 & 5 & 31.000 & 30.972 & 30.689 & 30.450 & 30.066 & 29.227 & 28.301 & 22.770 & 9.264 \\
 & 7 & 31.000 & 30.989 & 30.788 & 30.566 & 30.000 & 29.241 & 27.549 & 19.674 & 5.065 \\
 & 10 & 31.000 & 30.996 & 30.847 & 30.624 & 30.000 & 29.029 & 26.860 & 14.612 & 3.510 \\
 & 15 & 31.000 & 30.999 & 30.903 & 30.693 & 30.000 & 28.259 & 25.787 & 11.128 & 2.421 \\
\addlinespace[0.2em]
0.4 & 1 & 30.896 & 30.685 & 30.337 & 30.185 & 30.000 & 29.614 & 29.319 & 27.705 & 23.218 \\
 & 5 & 31.000 & 30.993 & 30.812 & 30.588 & 30.000 & 29.160 & 27.291 & 18.509 & 4.298 \\
 & 7 & 31.000 & 30.998 & 30.866 & 30.646 & 30.000 & 28.507 & 26.541 & 13.430 & 3.061 \\
 & 10 & 31.000 & 30.999 & 30.913 & 30.638 & 30.000 & 28.163 & 25.484 & 10.335 & 2.234 \\
 & 15 & 31.000 & 31.000 & 30.954 & 30.722 & 30.000 & 27.639 & 22.586 & 5.151 & 1.661 \\
\addlinespace[0.2em]
0.8 & 1 & 30.982 & 30.900 & 30.618 & 30.319 & 30.000 & 29.547 & 28.527 & 23.903 & 12.321 \\
 & 5 & 31.000 & 31.000 & 30.949 & 30.712 & 30.000 & 27.705 & 22.795 & 5.309 & 1.689 \\
 & 7 & 31.000 & 31.000 & 30.974 & 30.781 & 30.000 & 27.113 & 20.834 & 3.776 & 1.379 \\
 & 10 & 31.000 & 31.000 & 30.993 & 30.893 & 30.090 & 25.580 & 16.465 & 2.758 & 1.184 \\
 & 15 & 31.000 & 31.000 & 30.998 & 30.945 & 29.966 & 24.580 & 13.507 & 2.004 & 1.063 \\
\addlinespace[0.2em]
\end{longtable}
\endgroup

\begingroup
\setlength{\tabcolsep}{3pt}
\begin{longtable}{c c r r r r r r r r r}
\caption{$TARL_1$ values for the EWMA-RZ$^+$ chart with smoothing constant $\lambda = 0.1$, under unequal coefficients of variation ($\gamma_X \neq \gamma_Y$), horizon $I = 30$, target $TARL_0 = 30$, $z_0 = 1$, across nine shift levels $\tau$.}\label{tab:tarl1_unequal_T30_lam01}\\
\toprule
$\rho_0$ & $n$ & $\tau=0.9$ & $\tau=0.95$ & $\tau=0.98$ & $\tau=0.99$ & $\tau=1$ & $\tau=1.01$ & $\tau=1.02$ & $\tau=1.05$ & $\tau=1.1$ \\
\midrule
\endfirsthead
\multicolumn{11}{l}{\textit{Table \ref{tab:tarl1_unequal_T30_lam01} (continued)}}\\
\toprule
$\rho_0$ & $n$ & $\tau=0.9$ & $\tau=0.95$ & $\tau=0.98$ & $\tau=0.99$ & $\tau=1$ & $\tau=1.01$ & $\tau=1.02$ & $\tau=1.05$ & $\tau=1.1$ \\
\midrule
\endhead
\bottomrule
\endlastfoot
\multicolumn{11}{l}{\textit{$\gamma_X = 0.01$, $\gamma_Y = 0.2$}} \\
-0.8 & 1 & 30.938 & 30.769 & 30.467 & 30.337 & 30.000 & 29.769 & 29.489 & 27.780 & 21.563 \\
 & 5 & 31.000 & 30.999 & 30.925 & 30.706 & 30.000 & 29.069 & 26.219 & 14.791 & 3.557 \\
 & 7 & 31.000 & 31.000 & 30.922 & 30.755 & 30.000 & 27.744 & 25.277 & 7.793 & 2.715 \\
 & 10 & 31.000 & 31.000 & 30.954 & 30.805 & 30.000 & 27.298 & 21.081 & 5.868 & 2.104 \\
 & 15 & 31.000 & 31.000 & 30.979 & 30.859 & 30.000 & 26.634 & 19.016 & 4.226 & 1.627 \\
\addlinespace[0.2em]
-0.4 & 1 & 30.939 & 30.773 & 30.470 & 30.340 & 30.000 & 29.765 & 29.479 & 27.732 & 21.376 \\
 & 5 & 31.000 & 30.999 & 30.928 & 30.710 & 30.000 & 29.047 & 26.140 & 9.751 & 3.452 \\
 & 7 & 31.000 & 31.000 & 30.925 & 30.760 & 30.000 & 27.707 & 25.170 & 7.590 & 2.641 \\
 & 10 & 31.000 & 31.000 & 30.957 & 30.810 & 30.000 & 27.248 & 20.920 & 5.698 & 2.051 \\
 & 15 & 31.000 & 31.000 & 30.980 & 30.864 & 30.000 & 26.566 & 18.810 & 4.102 & 1.591 \\
\addlinespace[0.2em]
0.0 & 1 & 30.941 & 30.777 & 30.474 & 30.342 & 30.000 & 29.760 & 29.468 & 27.681 & 21.180 \\
 & 5 & 31.000 & 30.999 & 30.931 & 30.715 & 30.000 & 29.024 & 26.055 & 9.525 & 3.348 \\
 & 7 & 31.000 & 31.000 & 30.956 & 30.765 & 30.000 & 27.667 & 25.056 & 7.380 & 2.567 \\
 & 10 & 31.000 & 31.000 & 30.959 & 30.815 & 30.000 & 27.195 & 20.748 & 5.526 & 1.998 \\
 & 15 & 31.000 & 31.000 & 30.982 & 30.868 & 30.000 & 26.492 & 18.592 & 3.978 & 1.555 \\
\addlinespace[0.2em]
0.4 & 1 & 30.942 & 30.781 & 30.478 & 30.345 & 30.000 & 29.756 & 29.457 & 27.627 & 20.972 \\
 & 5 & 31.000 & 30.999 & 30.934 & 30.720 & 30.000 & 29.000 & 25.964 & 9.291 & 3.243 \\
 & 7 & 31.000 & 31.000 & 30.958 & 30.770 & 30.000 & 27.625 & 24.933 & 7.166 & 2.492 \\
 & 10 & 31.000 & 31.000 & 30.962 & 30.820 & 30.000 & 27.138 & 20.565 & 5.351 & 1.945 \\
 & 15 & 31.000 & 31.000 & 30.991 & 30.873 & 30.000 & 26.413 & 18.360 & 3.853 & 1.519 \\
\addlinespace[0.2em]
0.8 & 1 & 30.944 & 30.785 & 30.482 & 30.347 & 30.000 & 29.751 & 29.446 & 26.850 & 20.752 \\
 & 5 & 31.000 & 31.000 & 30.936 & 30.725 & 30.000 & 28.974 & 25.866 & 9.048 & 3.138 \\
 & 7 & 31.000 & 31.000 & 30.961 & 30.775 & 30.000 & 27.580 & 24.802 & 6.945 & 2.417 \\
 & 10 & 31.000 & 31.000 & 30.979 & 30.825 & 30.000 & 27.077 & 20.370 & 5.174 & 1.891 \\
 & 15 & 31.000 & 31.000 & 30.991 & 30.877 & 30.000 & 26.329 & 18.114 & 3.726 & 1.483 \\
\addlinespace[0.2em]
\midrule
\multicolumn{11}{l}{\textit{$\gamma_X = 0.2$, $\gamma_Y = 0.01$}} \\
-0.8 & 1 & 31.000 & 30.954 & 30.763 & 30.315 & 30.000 & 29.584 & 29.053 & 24.297 & 12.088 \\
 & 5 & 31.000 & 31.000 & 30.922 & 30.561 & 30.000 & 28.978 & 25.138 & 10.895 & 4.010 \\
 & 7 & 31.000 & 31.000 & 30.948 & 30.808 & 30.000 & 27.007 & 24.210 & 8.848 & 3.077 \\
 & 10 & 31.000 & 31.000 & 30.969 & 30.847 & 30.000 & 26.545 & 22.944 & 6.824 & 2.352 \\
 & 15 & 31.000 & 31.000 & 30.993 & 30.889 & 30.000 & 25.854 & 15.247 & 4.923 & 1.776 \\
\addlinespace[0.2em]
-0.4 & 1 & 31.000 & 30.956 & 30.552 & 30.320 & 30.000 & 29.576 & 29.032 & 24.194 & 11.896 \\
 & 5 & 31.000 & 31.000 & 30.925 & 30.568 & 30.000 & 28.953 & 25.045 & 10.671 & 3.899 \\
 & 7 & 31.000 & 31.000 & 30.950 & 30.812 & 30.000 & 26.964 & 24.091 & 8.629 & 2.994 \\
 & 10 & 31.000 & 31.000 & 30.971 & 30.851 & 30.000 & 26.489 & 22.791 & 6.632 & 2.292 \\
 & 15 & 31.000 & 31.000 & 30.994 & 30.892 & 30.000 & 25.778 & 15.050 & 4.777 & 1.735 \\
\addlinespace[0.2em]
0.0 & 1 & 31.000 & 30.958 & 30.560 & 30.325 & 30.000 & 29.566 & 29.009 & 24.086 & 11.697 \\
 & 5 & 31.000 & 31.000 & 30.929 & 30.576 & 30.000 & 27.281 & 24.946 & 10.438 & 3.787 \\
 & 7 & 31.000 & 31.000 & 30.953 & 30.816 & 30.000 & 26.919 & 23.964 & 8.405 & 2.910 \\
 & 10 & 31.000 & 31.000 & 30.973 & 30.855 & 30.000 & 26.430 & 22.627 & 6.437 & 2.231 \\
 & 15 & 31.000 & 31.000 & 30.994 & 30.896 & 30.000 & 25.698 & 14.843 & 4.629 & 1.694 \\
\addlinespace[0.2em]
0.4 & 1 & 31.000 & 30.961 & 30.568 & 30.331 & 30.000 & 29.557 & 27.371 & 23.970 & 11.492 \\
 & 5 & 31.000 & 31.000 & 30.932 & 30.584 & 30.000 & 27.243 & 24.840 & 10.196 & 3.675 \\
 & 7 & 31.000 & 31.000 & 30.956 & 30.821 & 30.000 & 26.870 & 23.828 & 8.174 & 2.826 \\
 & 10 & 31.000 & 31.000 & 30.975 & 30.859 & 30.000 & 26.366 & 22.452 & 6.239 & 2.170 \\
 & 15 & 31.000 & 31.000 & 30.995 & 30.900 & 30.000 & 25.612 & 14.624 & 4.481 & 1.652 \\
\addlinespace[0.2em]
0.8 & 1 & 31.000 & 30.963 & 30.576 & 30.337 & 30.000 & 29.546 & 27.340 & 23.846 & 11.278 \\
 & 5 & 31.000 & 31.000 & 30.935 & 30.592 & 30.000 & 27.202 & 24.726 & 9.944 & 3.561 \\
 & 7 & 31.000 & 31.000 & 30.958 & 30.655 & 30.000 & 26.817 & 23.682 & 7.936 & 2.741 \\
 & 10 & 31.000 & 31.000 & 30.977 & 30.864 & 30.000 & 26.297 & 22.263 & 6.037 & 2.108 \\
 & 15 & 31.000 & 31.000 & 30.995 & 30.904 & 30.000 & 25.519 & 14.392 & 4.332 & 1.611 \\
\addlinespace[0.2em]
\end{longtable}
\endgroup

\begingroup
\setlength{\tabcolsep}{3pt}
\begin{longtable}{c c r r r r r r r r r}
\caption{$TARL_1$ values for the EWMA-RZ$^+$ chart with smoothing constant $\lambda = 0.2$, under unequal coefficients of variation ($\gamma_X \neq \gamma_Y$), horizon $I = 30$, target $TARL_0 = 30$, $z_0 = 1$, across nine shift levels $\tau$.}\label{tab:tarl1_unequal_T30_lam02}\\
\toprule
$\rho_0$ & $n$ & $\tau=0.9$ & $\tau=0.95$ & $\tau=0.98$ & $\tau=0.99$ & $\tau=1$ & $\tau=1.01$ & $\tau=1.02$ & $\tau=1.05$ & $\tau=1.1$ \\
\midrule
\endfirsthead
\multicolumn{11}{l}{\textit{Table \ref{tab:tarl1_unequal_T30_lam02} (continued)}}\\
\toprule
$\rho_0$ & $n$ & $\tau=0.9$ & $\tau=0.95$ & $\tau=0.98$ & $\tau=0.99$ & $\tau=1$ & $\tau=1.01$ & $\tau=1.02$ & $\tau=1.05$ & $\tau=1.1$ \\
\midrule
\endhead
\bottomrule
\endlastfoot
\multicolumn{11}{l}{\textit{$\gamma_X = 0.01$, $\gamma_Y = 0.2$}} \\
-0.8 & 1 & 30.826 & 30.592 & 30.312 & 30.197 & 30.000 & 29.835 & 29.643 & 28.577 & 25.942 \\
 & 5 & 31.000 & 30.987 & 30.815 & 30.555 & 30.000 & 29.211 & 27.562 & 18.111 & 4.223 \\
 & 7 & 31.000 & 30.996 & 30.848 & 30.620 & 30.000 & 28.646 & 26.781 & 14.878 & 2.899 \\
 & 10 & 31.000 & 30.999 & 30.904 & 30.632 & 30.000 & 28.291 & 24.893 & 9.602 & 2.083 \\
 & 15 & 31.000 & 31.000 & 30.961 & 30.775 & 30.000 & 27.629 & 23.784 & 5.025 & 1.581 \\
\addlinespace[0.2em]
-0.4 & 1 & 30.829 & 30.596 & 30.315 & 30.199 & 30.000 & 29.832 & 29.637 & 28.551 & 25.842 \\
 & 5 & 31.000 & 30.988 & 30.819 & 30.560 & 30.000 & 29.193 & 27.505 & 17.839 & 4.071 \\
 & 7 & 31.000 & 30.996 & 30.853 & 30.626 & 30.000 & 28.620 & 26.699 & 14.564 & 2.803 \\
 & 10 & 31.000 & 30.999 & 30.908 & 30.639 & 30.000 & 28.254 & 24.767 & 9.318 & 2.024 \\
 & 15 & 31.000 & 31.000 & 30.963 & 30.781 & 30.000 & 27.575 & 23.610 & 4.843 & 1.544 \\
\addlinespace[0.2em]
0.0 & 1 & 30.832 & 30.600 & 30.318 & 30.201 & 30.000 & 29.830 & 29.534 & 28.523 & 25.736 \\
 & 5 & 31.000 & 30.989 & 30.824 & 30.566 & 30.000 & 29.175 & 27.444 & 17.553 & 3.918 \\
 & 7 & 31.000 & 30.997 & 30.858 & 30.632 & 30.000 & 28.592 & 26.612 & 12.255 & 2.708 \\
 & 10 & 31.000 & 30.999 & 30.912 & 30.646 & 30.000 & 28.214 & 24.632 & 9.027 & 1.966 \\
 & 15 & 31.000 & 31.000 & 30.965 & 30.787 & 30.000 & 27.517 & 22.125 & 4.661 & 1.508 \\
\addlinespace[0.2em]
0.4 & 1 & 30.835 & 30.605 & 30.321 & 30.203 & 30.000 & 29.827 & 29.527 & 28.493 & 25.623 \\
 & 5 & 31.000 & 30.990 & 30.829 & 30.571 & 30.000 & 29.155 & 27.379 & 17.250 & 3.767 \\
 & 7 & 31.000 & 30.997 & 30.863 & 30.638 & 30.000 & 28.561 & 25.799 & 11.932 & 2.614 \\
 & 10 & 31.000 & 30.999 & 30.916 & 30.654 & 30.000 & 28.172 & 24.486 & 6.489 & 1.908 \\
 & 15 & 31.000 & 31.000 & 30.968 & 30.793 & 30.000 & 27.455 & 21.914 & 4.479 & 1.472 \\
\addlinespace[0.2em]
0.8 & 1 & 30.839 & 30.609 & 30.325 & 30.205 & 30.000 & 29.824 & 29.520 & 28.462 & 25.504 \\
 & 5 & 31.000 & 30.991 & 30.834 & 30.577 & 30.000 & 29.134 & 27.309 & 16.930 & 3.616 \\
 & 7 & 31.000 & 30.997 & 30.868 & 30.645 & 30.000 & 28.529 & 25.689 & 11.595 & 2.521 \\
 & 10 & 31.000 & 31.000 & 30.920 & 30.661 & 30.000 & 28.127 & 24.330 & 6.244 & 1.850 \\
 & 15 & 31.000 & 31.000 & 30.970 & 30.799 & 30.000 & 27.389 & 21.687 & 4.297 & 1.437 \\
\addlinespace[0.2em]
\midrule
\multicolumn{11}{l}{\textit{$\gamma_X = 0.2$, $\gamma_Y = 0.01$}} \\
-0.8 & 1 & 30.999 & 30.944 & 30.662 & 30.473 & 30.000 & 29.570 & 29.013 & 25.480 & 16.233 \\
 & 5 & 31.000 & 30.998 & 30.876 & 30.666 & 30.000 & 28.962 & 26.396 & 15.061 & 3.606 \\
 & 7 & 31.000 & 31.000 & 30.913 & 30.715 & 30.000 & 28.722 & 25.550 & 12.355 & 2.700 \\
 & 10 & 31.000 & 31.000 & 30.958 & 30.767 & 30.000 & 28.389 & 24.360 & 6.760 & 2.048 \\
 & 15 & 31.000 & 31.000 & 30.979 & 30.826 & 30.000 & 27.076 & 22.528 & 4.645 & 1.569 \\
\addlinespace[0.2em]
-0.4 & 1 & 30.999 & 30.932 & 30.668 & 30.476 & 30.000 & 29.558 & 28.983 & 25.348 & 15.912 \\
 & 5 & 31.000 & 30.999 & 30.881 & 30.672 & 30.000 & 28.934 & 26.301 & 14.734 & 3.479 \\
 & 7 & 31.000 & 31.000 & 30.918 & 30.721 & 30.000 & 28.686 & 25.427 & 12.024 & 2.614 \\
 & 10 & 31.000 & 31.000 & 30.961 & 30.774 & 30.000 & 28.343 & 24.199 & 6.519 & 1.992 \\
 & 15 & 31.000 & 31.000 & 30.981 & 30.831 & 30.000 & 27.005 & 22.311 & 4.473 & 1.535 \\
\addlinespace[0.2em]
0.0 & 1 & 30.999 & 30.936 & 30.674 & 30.334 & 30.000 & 29.546 & 28.951 & 25.207 & 15.578 \\
 & 5 & 31.000 & 30.999 & 30.887 & 30.679 & 30.000 & 28.904 & 26.199 & 14.394 & 3.353 \\
 & 7 & 31.000 & 31.000 & 30.922 & 30.728 & 30.000 & 28.648 & 25.296 & 8.484 & 2.530 \\
 & 10 & 31.000 & 31.000 & 30.964 & 30.780 & 30.000 & 28.293 & 24.027 & 6.275 & 1.936 \\
 & 15 & 31.000 & 31.000 & 30.983 & 30.837 & 30.000 & 26.929 & 22.080 & 4.303 & 1.501 \\
\addlinespace[0.2em]
0.4 & 1 & 30.999 & 30.941 & 30.681 & 30.341 & 30.000 & 29.533 & 28.918 & 25.056 & 15.231 \\
 & 5 & 31.000 & 30.999 & 30.892 & 30.685 & 30.000 & 28.871 & 26.089 & 14.038 & 3.228 \\
 & 7 & 31.000 & 31.000 & 30.927 & 30.734 & 30.000 & 28.607 & 25.154 & 8.186 & 2.446 \\
 & 10 & 31.000 & 31.000 & 30.956 & 30.787 & 30.000 & 27.503 & 23.841 & 6.030 & 1.882 \\
 & 15 & 31.000 & 31.000 & 30.984 & 30.843 & 30.000 & 26.847 & 21.832 & 4.133 & 1.467 \\
\addlinespace[0.2em]
0.8 & 1 & 30.999 & 30.945 & 30.592 & 30.349 & 30.000 & 29.519 & 28.881 & 24.894 & 14.869 \\
 & 5 & 31.000 & 30.999 & 30.898 & 30.692 & 30.000 & 28.836 & 25.970 & 13.665 & 3.106 \\
 & 7 & 31.000 & 31.000 & 30.931 & 30.741 & 30.000 & 28.562 & 25.000 & 7.881 & 2.363 \\
 & 10 & 31.000 & 31.000 & 30.959 & 30.793 & 30.000 & 27.438 & 23.641 & 5.782 & 1.827 \\
 & 15 & 31.000 & 31.000 & 30.986 & 30.850 & 30.000 & 26.759 & 21.567 & 3.964 & 1.435 \\
\addlinespace[0.2em]
\end{longtable}
\endgroup

\fig{fig:profile} provides a visual summary of the $\mathrm{TARL}_1$ profile for one slice of the design space. The five sample sizes give curves of the same shape, all peaking around $\tau=1$ at $\mathrm{TARL}_1\approx\mathit{TARL}_0^*$ and decreasing toward $1$ as $\tau$ grows. The vertical spread of the curves quantifies the benefit of larger samples: at $\tau=1.05$, increasing $n$ from $1$ to $15$ divides $\mathrm{TARL}_1$ by a factor of approximately three.

\subsection{Comparison with the ShRZ chart of \citet{Tran2021}}
\label{ssec:comparison}

For a like-for-like comparison with the Shewhart-type ShRZ chart of \citet{Tran2021}, we use the same parameter grid and the same target $\mathit{TARL}_0^* = I$. Across the upward-shift levels $\tau\in\{1.01,1.02,1.05\}$, which are arguably the most relevant in practice, the proposed EWMA-RZ$^+$ chart reduces $\mathrm{TARL}_1$ by approximately 15--40\% compared to the ShRZ chart in the equal-CV cases with moderate-to-large CVs. The two charts are essentially equivalent in the large-shift regime ($\tau\ge 1.10$), where both signal within one or two inspections. The improvement at small shifts is the well-known consequence of the memory effect of the EWMA statistic: small persistent shifts accumulate in $W_i$ and ultimately push it above $\mathit{UCL}$, while a Shewhart chart, which depends only on the current $\hat Z_i$, remains insensitive until the shift is large enough that a single observation crosses the limit on its own.

%%%%%%%%%%%%%%%%%%%%%%%%%%%%%%%%%%%%%%%%%%%%%%%%%%%%%%%%%%%%%%%%%%%%%%%%%%%%%%%%
%%  >>>>> INSERTION 2: new subsection 5.6 <<<<<
%%%%%%%%%%%%%%%%%%%%%%%%%%%%%%%%%%%%%%%%%%%%%%%%%%%%%%%%%%%%%%%%%%%%%%%%%%%%%%%%
\subsection{Robustness to departures from bivariate normality}
\label{ssec:robustness}

The calibration of the upper control limit in Section~\ref{sec:optimization}, and the entire factorial study of Sections~\ref{ssec:limits}--\ref{ssec:comparison}, rest on the bivariate normal model of equation~\eqref{eq:bvn}. To assess how sensitive the chart's actual performance is to this modelling choice, we re-evaluate $\mathrm{TARL}_0$ and $\mathrm{TARL}_1$ at the \emph{same} calibrated $\mathit{UCL}$ under three families of non-normal data-generating processes, each matched to the in-control bivariate normal in mean, coefficient of variation, and Pearson correlation $\rho_0$:
\begin{itemize}
	\item bivariate lognormal, obtained by exponentiating a bivariate normal whose underlying correlation is solved so that the Pearson correlation of $(X,Y)$ equals $\rho_0$;
	\item bivariate Student-$t$ with $\nu=10$ degrees of freedom (moderate tail heaviness);
	\item bivariate Student-$t$ with $\nu=5$ degrees of freedom (heavy tails, excess kurtosis $6$).
\end{itemize}
For each pair $(\tau,\text{distribution})$ we generate $N_{\mathrm{MC}}=5\times 10^{5}$ independent short runs of $I$ inspections each, apply the EWMA recursion~\eqref{eq:EWMA} with the original $\mathit{UCL}$, and average the run lengths. The bivariate-normal baseline is included for cross-comparison with the Markov-chain values of Sections~\ref{ssec:limits}--\ref{ssec:detection} and uses the same sample sizes throughout; minor differences ($\le 3\%$) between the direct-MC baseline and the Markov-chain value reported earlier reflect the closed-form approximation~\eqref{eq:Zcdf}, which we have already noted is used here as a fast surrogate.

Tables~\ref{tab:robustness_A} and~\ref{tab:robustness_B} report the results for two representative configurations chosen to span the CV range used in the paper. Configuration~A mirrors the low-CV regime of the illustrative example of Section~\ref{sec:illustrative} ($\gamma_X=\gamma_Y=0.05$, $I=20$, $\mathit{UCL}=1.01918$). Configuration~B uses the largest CV pair of the factorial study ($\gamma_X=\gamma_Y=0.20$, $I=10$, $\mathit{UCL}=1.0621$), where any sensitivity to non-normality should be most pronounced. In both regimes, $\lambda=0.2$, $\rho_0=0.4$, $n=5$, $z_0=1$, and the target in-control $\mathit{TARL}_0$ equals $I$.

\begin{table}[!htbp]
	\centering
	\caption{$\mathrm{TARL}$ values of the EWMA-RZ$^+$ chart under four data-generating processes with matched first two moments and matched correlation. Low-CV configuration: $\gamma_X=\gamma_Y=0.05$, $\rho_0=0.4$, $\lambda=0.2$, $n=5$, $I=20$, $\mathit{UCL}=1.01918$, $z_0=1$. Monte Carlo with $5\times 10^{5}$ replications per cell; maximum standard error $0.009$.}
	\label{tab:robustness_A}
	\begin{tabular}{lcccccc}
		\toprule
		Distribution & $\tau=0.95$ & $\tau=1.00$ & $\tau=1.01$ & $\tau=1.02$ & $\tau=1.05$ & $\tau=1.10$ \\
		\midrule
		Normal (baseline)         & 21.000 & 20.087 & 15.462 & 8.772 & 2.837 & 1.445 \\
		Lognormal                 & 21.000 & 20.080 & 15.465 & 8.777 & 2.837 & 1.446 \\
		Student-$t$ ($\nu=10$)    & 21.000 & 20.042 & 15.475 & 8.798 & 2.835 & 1.442 \\
		Student-$t$ ($\nu=5$)     & 20.998 & 19.975 & 15.582 & 8.883 & 2.831 & 1.442 \\
		\bottomrule
	\end{tabular}
\end{table}

\begin{table}[!htbp]
	\centering
	\caption{$\mathrm{TARL}$ values of the EWMA-RZ$^+$ chart under four data-generating processes with matched first two moments and matched correlation. High-CV configuration: $\gamma_X=\gamma_Y=0.20$, $\rho_0=0.4$, $\lambda=0.2$, $n=5$, $I=10$, $\mathit{UCL}=1.0621$, $z_0=1$. Monte Carlo with $5\times 10^{5}$ replications per cell; maximum standard error $0.005$.}
	\label{tab:robustness_B}
	\begin{tabular}{lcccccc}
		\toprule
		Distribution & $\tau=0.95$ & $\tau=1.00$ & $\tau=1.01$ & $\tau=1.02$ & $\tau=1.05$ & $\tau=1.10$ \\
		\midrule
		Normal (baseline)         & 10.929 & 10.206 & 9.844 & 9.400 & 7.604 & 4.670 \\
		Lognormal                 & 10.945 & 10.256 & 9.898 & 9.453 & 7.655 & 4.673 \\
		Student-$t$ ($\nu=10$)    & 10.915 & 10.190 & 9.844 & 9.403 & 7.639 & 4.680 \\
		Student-$t$ ($\nu=5$)     & 10.876 & 10.188 & 9.860 & 9.428 & 7.691 & 4.711 \\
		\bottomrule
	\end{tabular}
\end{table}

Two observations stand out. First, in the low-CV regime (\tab{tab:robustness_A}), all four distributions give essentially indistinguishable $\mathrm{TARL}$ profiles: the largest cell-wise discrepancy is below $0.12$, well within the $\pm 0.01$ Monte Carlo precision band. This is consistent with the central limit theorem: with $n=5$ and small $\gamma$, the sample means $\bar X_i$ and $\bar Y_i$ are very close to normal regardless of the marginal distribution of the individual observations, so the ratio $\hat Z_i$ also behaves nearly identically across the four distributions. Second, in the high-CV regime (\tab{tab:robustness_B}), the four profiles remain remarkably similar. Lognormal data yields a slightly \emph{larger} $\mathrm{TARL}_0$ ($+0.5\%$) than the normal baseline, meaning a marginally lower false-alarm rate; the heavy-tailed Student-$t$ with $\nu=5$ yields a marginally smaller $\mathrm{TARL}_0$ ($-0.2\%$). Out-of-control performance differs by at most $1\%$ across the four distributions at any shift level. In operational terms, the chart calibrated under the bivariate normal model preserves both its prescribed false-alarm rate and its detection power when the data exhibit moderate departures from normality of the kinds tested.

We caution that the robustness demonstrated here is for the moderate non-normality scenarios listed above, with matched first two moments and matched correlation, and with $n=5$ providing some central-limit smoothing of the sample means. Stronger departures (e.g.\ Student-$t$ with $\nu\le 4$, or pronounced marginal skewness combined with $n=1$) or asymmetric tail dependence captured by non-Gaussian copulas may warrant separate treatment; we list this as a direction for future work in Section~\ref{sec:conclusions}.

%%%%%%%%%%%%%%%%%%%%%%%%%%%%%%%%%%%%%%%%%%%%%%%%%%%%%%%%%%%%%%%%%%%%%%%%%%%%%%%%
\section{Illustrative example}
\label{sec:illustrative}

To illustrate the practical use of the EWMA-RZ$^+$ chart, consider a beverage filling line that bottles a carbonated soft drink. Each bottle is characterized by two quality variables: $X$, the dissolved carbon-dioxide concentration (g/L), and $Y$, the filled volume (mL). Both variables are routinely controlled on the line, but it is the ratio $Z=X/Y$, i.e.\ the carbon-dioxide content per unit of volume, that determines the consumer-perceived fizziness of the product and that is therefore the relevant quality indicator for marketing specifications. The line is operated under SPR conditions because the same equipment is used for several recipes, with each recipe being run for a short batch before changeover. Suppose that for the current recipe the in-control specifications are $z_0=1$ (in normalized units), $\gamma_X=\gamma_Y=0.05$, $\rho_0=0.4$. The production run consists of $I=20$ inspections, with $n=5$ bottles sampled at each inspection. The smoothing constant is set to $\lambda=0.2$, and the target in-control TARL is $\mathit{TARL}_0^*=20$.

\begin{figure}[H]
	\centering
	\includegraphics[width=0.85\linewidth]{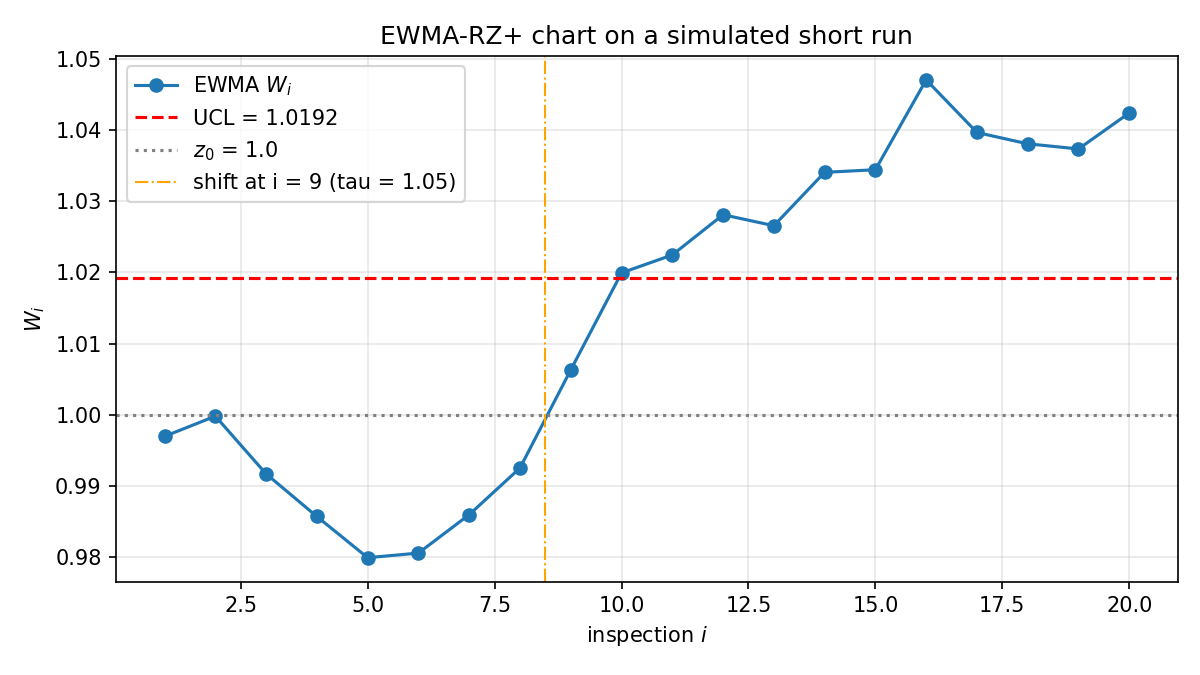}
	\caption{Illustrative EWMA-RZ$^+$ chart for the beverage-filling example. The process drifts upward to $\tau=1.05$ at inspection $i=9$ and the chart signals at $i=10$. Calibrated $\mathit{UCL}=1.01918$, $z_0=1$, $\gamma_X=\gamma_Y=0.05$, $\rho_0=0.4$, $\lambda=0.2$, $n=5$, $I=20$, $\mathit{TARL}_0^*=20$.}
	\label{fig:illustrative}
\end{figure}

Applying the design procedure of Section~\ref{sec:optimization} yields a calibrated upper control limit
\begin{align*}
	\mathit{UCL} = 1.01918,
	\qquad
	\mathrm{TARL}_0(\mathit{UCL}) = 20.00.
\end{align*}
We then simulate a production run on this recipe. From inspection $i=1$ to $i=8$, the process operates in control with $\mu_X/\mu_Y=z_0=1$. At inspection $i=9$ a mechanical drift in the carbonator induces an upward shift to $\tau=1.05$ in the ratio. The resulting EWMA chart is shown in \fig{fig:illustrative}. After a few inspections of fluctuation around $z_0$, the EWMA statistic crosses the calibrated $\mathit{UCL}$ at inspection $i=10$, i.e.\ within two samples of the actual onset of the shift. This rapid detection allows the operator to interrupt the run, recalibrate the carbonator, and resume production with a minimal volume of out-of-specification product.

For comparison, the ShRZ Shewhart-type chart of \citet{Tran2021} applied to the same simulated run signals only later: a single $\hat Z_i$ at the post-shift level $\tau z_0 = 1.05$ rarely exceeds the Shewhart limit, so the memory effect of the EWMA is decisive in keeping the time-to-signal short.

%%%%%%%%%%%%%%%%%%%%%%%%%%%%%%%%%%%%%%%%%%%%%%%%%%%%%%%%%%%%%%%%%%%%%%%%%%%%%%%%
\section{Conclusions}
\label{sec:conclusions}

This paper proposed an EWMA-type control chart for monitoring the ratio of two normally distributed quality characteristics under short production run conditions. The construction relies on (i) the corrected closed-form density of the ratio derived by \citet{Nadarajah2020}, with the closed-form approximation of \citet{Celano2016_Synthentic_RZ} used as a fast surrogate for the Markov-chain build; (ii) an upper-sided EWMA recursion on the ratio of sample means; and (iii) a vectorized Markov-chain representation that turns the upper-control-limit calibration into a one-dimensional root-finding problem solvable in milliseconds.

An extensive factorial study covering the smoothing constant $\lambda$, the in-control correlation $\rho_0$, the coefficients of variation $(\gamma_X,\gamma_Y)$, the sample size $n$, the horizon $I$, and the shift factor $\tau$ shows that the proposed EWMA-RZ$^+$ chart calibrates accurately to the prescribed in-control TARL and detects small to moderate upward shifts substantially faster than the Shewhart-type ShRZ chart of \citet{Tran2021}. The improvement is most pronounced for shifts in the range $\tau\in[1.01,1.05]$, which is the most relevant regime for early intervention in industrial practice. An illustrative beverage-filling example demonstrates that the chart signals within two inspections of the onset of a $5\%$ upward drift, which would have been missed by a Shewhart-type design. A complementary robustness study (Section~\ref{ssec:robustness}) further indicates that the chart calibrated under bivariate normality preserves both its false-alarm rate and its detection power when the data follow lognormal or moderately heavy-tailed Student-$t$ distributions matched in mean, CV, and correlation.

Several extensions are natural for future work. First, the upper-sided design can be combined with a symmetric lower-sided design to yield a two-sided EWMA-RZ chart, with the TARL computed from the joint Markov chain of the two statistics. Second, adaptive sampling strategies such as variable sampling intervals (VSI) or variable sample size (VSS), already explored for related charts \citep{Castagliola2013,Amdouni2015,Nguyen2019}, would further enhance the chart's responsiveness. Third, a fully economic-statistical design of the chart, balancing sampling cost and false-alarm cost over the finite horizon, could be developed along the lines of \citet{Celano2012,Zhang2014}. Fourth, replacing the closed-form approximation \eqref{eq:Zcdf} by direct numerical integration of the exact density \eqref{eq:fz} in the Markov-chain build would allow a tightening of the accuracy of the calibration in the regime of small $n$ and large CVs, at a modest computational cost.
%%%%%%%%%%%%%%%%%%%%%%%%%%%%%%%%%%%%%%%%%%%%%%%%%%%%%%%%%%%%%%%%%%%%%%%%%%%%%%%%
%%  >>>>> INSERTION 3: new ``Fifth'' future-work direction <<<<<
%%%%%%%%%%%%%%%%%%%%%%%%%%%%%%%%%%%%%%%%%%%%%%%%%%%%%%%%%%%%%%%%%%%%%%%%%%%%%%%%
Fifth, a more comprehensive robustness analysis would complement the moderate-departure study of Section~\ref{ssec:robustness}. Three extensions are particularly relevant: (i) heavier-tailed marginal distributions such as Student-$t$ with very low degrees of freedom or contaminated-normal mixtures; (ii) bivariate copulas allowing asymmetric tail dependence between $X$ and $Y$, which the Pearson correlation does not fully capture; and (iii) more pronounced marginal skewness, especially at $n=1$ where central-limit smoothing of the sample means is absent. When the closed-form approximation~\eqref{eq:Zcdf} becomes inaccurate in any of these settings, $\mathit{UCL}$ can be re-calibrated by direct Monte Carlo on the assumed data-generating process at modest additional cost; the Markov-chain build of Section~\ref{sec:optimization} then plays the role of a fast initial estimate that the Monte Carlo refines.
%%%%%%%%%%%%%%%%%%%%%%%%%%%%%%%%%%%%%%%%%%%%%%%%%%%%%%%%%%%%%%%%%%%%%%%%%%%%%%%%
Finally, a robust version of the chart that accounts for estimation error in the in-control parameters $(\gamma_X,\gamma_Y,\rho_0)$ from a Phase~I sample of limited size would be an important contribution toward operational deployment in genuinely data-scarce SPR settings.

\section*{Data availability statement}

No external empirical dataset was used in this study. The results are based on analytical calculations and numerical simulations. The simulated data and Python code used to reproduce the findings are available from the corresponding author upon reasonable request.

\section*{Declaration of Competing Interest}

The authors declare that they have no known competing financial interests or personal relationships that could have influenced the work reported in this paper.

\bibliography{paper}
%%%%%%%%%%%%%%%%%%%%%%%%%%%%%%%%%%%%%%%%%%%%%%%%%%%%%%%%%%%%%%%%%%%%%%%%%%%%%%%%
\end{document}